\documentclass[10pt,reqno]{amsart}
\usepackage[backref=page]{hyperref}
\usepackage{amsaddr}
\usepackage{graphicx}
\usepackage{amssymb}
\usepackage{amsfonts}
\usepackage[]{amsmath}
\usepackage[]{epsfig}
\usepackage[]{pstricks}
\newpsobject{malla}{psgrid}{subgriddiv=1,griddots=10,gridlabels=6pt}
\usepackage[]{float}
\usepackage{setspace}
\usepackage{tikz}
\usepackage{listings}
\usepackage{comment}
\usepackage{subfigure}
\usepackage{orcidlink}

\usetikzlibrary{decorations.shapes}
\newpsobject{malla}{psgrid}{subgriddiv=1,griddots=10,gridlabels=6pt}
\newtheorem{theorem}{Theorem}[section]
\newtheorem{lemma}{Lemma}[section]

\newtheorem{proposition}{Proposition}[section]

\newtheorem{example}{Example}[section]
\numberwithin{equation}{section}
\theoremstyle{remark}
\newtheorem{remark}{Remark}[section]
\newcommand{\R}{\mathbb R}

\newcommand\comentario[1]\null
\begin{document}
\pagenumbering{arabic}	
\title[KdV on general metric star graphs]{The Korteweg-de Vries equation on general star graphs}
\author{M\'arcio CAVALCANTE \orcidlink{0000-0001-9873-1765}}
\address{\emph{Instituto de Matem\'atica, Universidade Federal de Alagoas,\\ Macei\'o-Brazil}}
\email{marcio.melo@im.ufal.br}
\let\thefootnote\relax\footnote{2020 Mathematics Subject Classification (AMS): 35Q53.}
\footnote{Permanent address: Instituto de Matem\'atica, Universidade Federal de Alagoas,
	Macei\'o-Brazil
}

\author{José Marques \orcidlink{0009-0001-9009-2628}}
\address{\emph{Instituto de Matem\'atica, Universidade Federal de Alagoas,\\ Macei\'o-Brazil}}
\email{jose.neto@im.ufal.br}

\begin{abstract}
In this paper, we establish local well-posedness for the Cauchy problem associated with the Korteweg-de Vries (KdV) equation on a general metric star graph. The graph comprises $m+k$ semi-infinite edges: $k$ negative half-lines and $m$ positive half-lines, all joined at a common vertex. The choice of boundary conditions is compatible with the conditions determined by the semigroup theory considered by Mugonolo, Noja, and Seifert (APDE 2018). The crucial point in this work is to obtain the integral formula using the forcing operator method and the Fourier restriction method of Bourgain (GAFA 1993). This study extends the results obtained by Cavalcante (ZAMP 2018) for the specific case of the $\mathcal Y$ junction to a more general class of star graphs.

\end{abstract}
\maketitle

\tableofcontents

\section{Introduction}

The well-known Korteweg-de Vries (KdV) equation
\begin{equation}\label{KdV}
u_t+u_{xxx}+u_xu=0,\ x,t\in \R
\end{equation}
was introduced by Boussinesq \cite{bous} in 1877 as a model for long waves propagating on a shallow water surface. Boussinesq's primary aim was to derive an equation with a smooth traveling wave solution, inspired by the wave observed in 1834 by Scott Russell \cite{SR}. 
 Later, it was rediscovered by Korteweg and de Vries \cite{KDV}.  The KdV equation is now commonly accepted as a mathematical model
for the unidirectional propagation of small-amplitude long waves in nonlinear dispersive
systems. In particular, the KdV equation is not only used to serve as a model to study surface
water waves. In fact, recently the KdV equation has been used as a model to study blood pressure waves in large arteries. In this way, Shi et al. \cite{Shi} highlight that while most 1D cardiovascular models neglect vessel wall inertia, some researchers (see \cite{Crepeau,Demiray,Laleg,Yomosa}) have included this factor, leading to a Korteweg-de Vries equation and solitary wave dynamics. They suggested that this approach may better represent arterial pulse waves and warrants further investigation.

In the mathematical context the Cauchy problem for the KdV posed on the
real axis, torus, on the half-lines and on a finite interval has been well studied in the last years, and we recommend the papers cited in \cite{MC} for a sample of this subject. 

The advance of studies concerning the KdV equation on star graphs is slow, if compared with other domains. The linearized version  was studied in the works \cite{Sobirov1,Sobirov2,Sobirov3}.
Also, Mugnolo, Noja and Seifert \cite{Noja} obtained a characterization of all boundary conditions under which the Airy-type evolution equation generates a contractive or unitary semigroup of operators. In \cite{MC}, Cavalcante obtained local well-posedness for the Cauchy problem associated with Korteweg–de Vries equation on a metric star graph with three semi-infinite edges given by one negative half-line and two positive half-lines attached to a common vertex, for two classes of boundary conditions. More recently, Parada, Cr\'epeau and Prieur \cite{Parada1} studied stabilization of KdV equation in a network with bounded and unbounded lengths. In \cite{solitons2} the authors studied, for Korteweg-de Vries equation on a metric star graph, existence of solitary waves in terms of the coefficients of the equation on each edge, the coupling condition at the central vertex of the star and the speeds of the travelling wave. Finally, results for the KdV equations on star graphs with bounded edges we cite the works \cite{Parada2,Parada3,Parada5}. Finally, linear and nonlinear instability for stationary solutions was obtained in \cite{Angulo1,Angulo2}.

 \subsection{Airy Equation on star graphs}\label{sub1}
In this section, to motivate the setting of our problem, we focus on the linear component of the Korteweg-de Vries equation, namely the Airy equation, within the context of metric star graphs featuring semi-infinite edges. This specific configuration has been previously explored in \cite{Noja}. The investigation of the Airy operator on star graphs serves as a crucial foundation for motivating our chosen problem. This is due to the fact that selecting appropriate boundary conditions for studying the KdV equation on star graphs presents a significant challenge.

Here, we consider a star graph $\mathcal{G},$ composed by $k+m$ edges in a set of edges $E=E_{-}\cup E_{+}$, where $E_{-}$ is a finite set of negative half-lines  and $E_{+}$ a finite set of positive half-lines.

Then, we consider the Airy equation 

$$\partial_t u_e = - \partial_x^3 u_e   \text{ on } \mathbb{R} \times \mathcal G.$$

We revisit some existence results for the Airy equation on star graphs, which were obtained by Mugnolo, Noja, and Seifert \cite{Noja} using semigroup theory. We also recommend the work \cite{solitons2}, where the same authors explored the existence of solitons for KdV equations on star graphs for some boundary conditions.

All vector spaces appearing are over the field $\mathbb{K} \in \{\mathbb{R}, \mathbb{C}\}$. Define the minimal operator $A_0$ in $H :=L^2(\mathcal G):= \bigoplus_{e \in E_-} L_2(-\infty, 0) \bigoplus_{e \in E_+} L_2(0,\infty)$ by 

$$D(A_0) := \bigoplus_{e \in E_-} C_c^\infty (-\infty, 0) \bigoplus_{e \in E_+} C_c^\infty (0,\infty),$$

$$A_0 u := (- \partial_x^3u_e )_{e \in E} \quad (u \in D(A_0)).$$
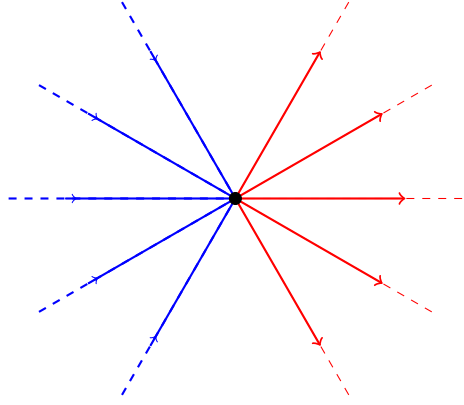
\begin{figure}
\begin{center}
\begin{tikzpicture}[scale=1.5]
\foreach \angle in {0,30,60,300,330} {
    \draw[dashed, red] (\angle:2cm) -- (0,0);
    \draw[->, thick, red] (\angle:0cm) -- (\angle:1.5cm);
}

\foreach \angle in {120,150,180,210,240} {
    \draw[dashed,thick,blue] (\angle:2cm) -- (0,0);
    \draw[->,blue] (\angle:1.5cm) -- (\angle:1.4cm);
    \draw[-,thick, blue] (\angle:1.5cm) -- (\angle:0.0cm);
}
\draw[fill] (0,0) circle (1.5pt) node[below right] {};
\end{tikzpicture}\caption{Star graph with five negative half-lines (the blue ones) and five positive half-lines (the red ones)}
\end{center}
\end{figure}

A simple computation proves that 

$$
D(A_0^*) = \bigoplus_{e \in E_-} H^3(-\infty, 0) \bigoplus_{e \in E_+} H^3(0,\infty),$$

$$A_0^* u = (-\partial_x^3u_e)_{e \in E} \quad (u \in D(A_0^*)),$$
and hence $A_0$ is skew-symmetric, i.e. densely defined and $A_0 \subset -A_0^*$.

An important point is finding operator $A$ such that $A_0 \subset A \subset -A_0^*$ and for which the abstract Cauchy problem
$$
\begin{cases}
\partial_t u_e(t) = Au_e(t),& (t > 0),\\
u_e(0) = u_0\in H
\end{cases}
$$
is well-posed, i.e. $A$ is the generator of a $C_0$-semigroup. To do this, we will have to impose coupling conditions at the vertex $0$. Note that this abstract Cauchy problem is a functional analytic description of the Airy equation at each edge $e \in E$, including the vertex conditions at $0.$

For $u \in D(A_0^*)$ let us abbreviate 

$$
\partial_x^ju(0-) := (\partial_x^ju_e(0-))_{e \in E_-}, \quad \partial_x^ju(0+) := (\partial_x^ju(0+))_{e \in E_+} \quad (j \in \{0,1,2\}).$$

Moreover, let us denote the identity operators on $\mathbb K^{k+m}$ defined by

$$
\mathbf{1}_{\pm} x := ( x_e)_{e \in E_{\pm}}.
$$

Define the spaces
$$
G_{-} := \mathbb K^{3k}\ \text{and}\ G_{+} := \mathbb K^{3m}
$$

and the block operator matrices $\mathcal B_{\pm}$ on $G_{\pm}$ by

\begin{equation}
\mathcal B_{\pm} = \begin{pmatrix}
0 &0 & \mathbf{1}_{\pm} \\
0 &-\mathbf{1}_{\pm} & 0 \\
\mathbf{1}_{\pm} & 0 & 0
\end{pmatrix}.
\end{equation}

We have $\mathcal B_{\pm} \in \mathcal{L}(G_{\pm})$, $\mathcal B_{\pm}$ is injective, and $\mathcal B_{\pm}^{-1} \in \mathcal{L}(G_{\pm})$. Define an indefinite inner product $\langle \cdot, \cdot \rangle_+ : G_+ \times G_+ \to \mathbb K$ by

$$
\langle x, y \rangle_{\pm} := (\mathcal B_{\pm} x | y)_{G_{\pm}} \quad (x,y \in G_{\pm}).
$$

Then $\mathbf K_{\pm} := (G_{\pm}, \langle \cdot, \cdot \rangle_{\pm})$ are Krein spaces.

For a linear operator $L: \mathbf K_- \to \mathbf K_+$, we define the following
\[
D(A_L) := \{u \in D(A); L(u(0-), \partial_xu(0-), \partial_x^2 u(0-)) = (u(0+), \partial_x u(0+), \partial_x ^2u(0+))\},
\]
\[
A_L u := -Au \quad (u \in D(A_L)).
\]

The work \cite{Noja} gives a complete characterization of the boundary conditions, which implies a nice theory for the Airy equation, depending on some properties of the operator $L$.

A special case of coupling conditions at the vertex $0$, obtained in \cite{Noja}, and described in \cite{solitons2}, stems from separating the boundary values of the first derivatives. In order to formulate this, let $Y \subset \mathbb K^{k} \oplus \mathbb K^m$ be a subspace and $B: \mathbb K^m  \to \mathbb K^k$ be linear. We define 
\[
\begin{split}D(A_{Y,B}) := \{u \in D(A_0^*);& (u(0-), u(0+)) \in Y,\\
& {\begin{pmatrix}
-\mathbf{1}_{-} & 0 \\
0 & \mathbf{1}_{+}
\end{pmatrix}
\begin{pmatrix}
u''(0-) \\
u''(0+)
\end{pmatrix}
\in Y^{\perp},}\\
& u'(0-) = Bu'(0+)\},\\
A_{Y,B} u :=& -A_0^*u \quad (u \in D(A_{Y,B})).
\end{split}
\]
\begin{proposition}\label{MNS}(Corolary 3.20 and Theorem 3.23 \cite{Noja} or Proposition 3.2 \cite{solitons2}). Let  $Y \subset \mathbb K^{k} \oplus \mathbb K^m$ be a subspace and $B: \mathbb K^{m} \to \mathbb K^k$ linear. Then:
\begin{itemize}
\item[(a)] $A_{Y,B}$ generates a unitary $C_0$ group if and only if $B$ is unitary.
\item [(b)] $A_{Y,B}$ generates a contractive $C_0$ semigroup if and only if $B$ is a contraction matrix. 
\end{itemize}
\end{proposition}

In the next subsection, we consider a specific class of boundary conditions, with derivatives of each order separated, that are included in the class of conditions given by Proposition \ref{MNS}.

\subsection{Formulation of the problem} 

The main purpose of this paper is to generalize the study of the Korteweg-de Vries equation considered in previous work \cite{MC} for a more general star graph, by considering a class of boundary conditions separating the derivatives.
 
 In this spirit, motivated by  \cite{Noja,Sobirov2}  we consider the KdV equation on a general star graph given by $k$ negative half-lines and $m$ positive half-lines attached to a common point, called the vertex. We denote this graph by $\mathcal G$ composed by $k$ negatives half-lines and $m$ positive half-lines, attached by an unique vertex. The vertex corresponds to $0$. 
 
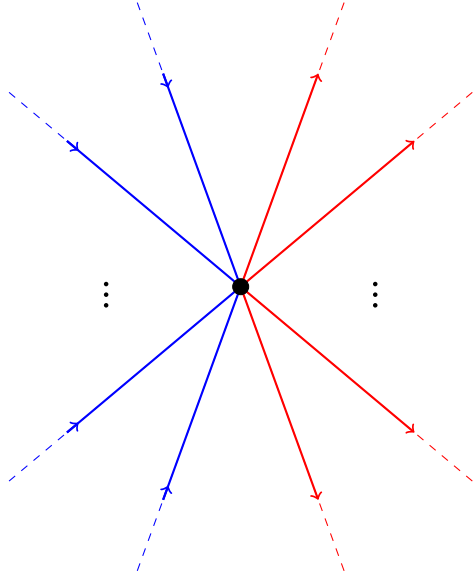
\begin{figure}[htp]\label{figure1}
	\centering 
	\begin{tikzpicture}[scale=2]
\node at (0.89,0.0)[rotate=0]{\Huge$\vdots$};
\node at (-0.89,0.0)[rotate=0]{\Huge$\vdots$};
	
\foreach \angle in {40,70,290,320} {
    \draw[dashed, red] (\angle:2cm) -- (0,0);
    \draw[->,thick,red] (\angle:0cm) -- (\angle:1.5cm);
}

\foreach \angle in {110,140,220,250} {
    \draw[dashed, blue] (\angle:2cm) -- (0,0);
    \draw[->,thick,blue] (\angle:1.5cm) -- (\angle:1.4cm);
    \draw[-,thick, blue] (\angle:1.5cm) -- (\angle:0.0cm);
}
\draw[fill] (0,0) circle (1.5pt) node[below right] {};
	\end{tikzpicture}
	\caption{A star graph with $m+k$ edges}
	\end{figure}

Based on Proposition \ref{MNS} of the previous discussion, we consider the following Cauchy problem associated to KdV equation on star graphs $\mathcal{G}$ given by
\begin{equation}\label{KDV}
\begin{cases}
\partial_t u_i+\partial^3_{x}u_i+u_i\partial_{x}u_i=0,& (x, t)\in (-\infty,0)\times (0,T),\ i=1,2, ..., k;\\
\partial_t v_j+\partial^3_{x}v_j+v_j\partial_{x}v_j=0,& (x,t)\in (0,\infty)\times (0,T),\ j=1,2, ..., m;\\
u_i(x,0)=u_{0,i}(x),&  x \in  (-\infty,0), i=1,2,..., k;\\
v_j(x,0)=v_{0,j}(x),&  x \in  (0,\infty), j=1,2,..., m;\\
\end{cases}
\end{equation}
with the following set of boundary conditions.
\begin{equation}\label{condfront}
\begin{cases}
u_1(0,t)=a_i u_i(0,t)=a_{k+j} v_j(0,t),\ i=2, ..., k,\ j=1, ..., m;& t\in (0,T)\\
\partial_x U(0,t)=B\partial_xV(0,t),& t\in (0,T)\\
\displaystyle\sum^k_{i=1} a_i^{-1}\partial^2_{x} u_i(0,t)=\displaystyle\sum^m_{j=1} a^{-1}_{k+j} \partial^2_{x} v_j(0,t) & t\in (0,T),
\end{cases}
\end{equation}
where   $V(x,t)=(v_1(x,t), ..., v_{m}(x,t))^T,\ U(x,t)=(u_1(x,t), ..., u_{k}(x,t))^T$ and $B=\left[b_{ij}\right]_{k\times m}$ is a $k \times m$ matrix such that  $$\| B(x)\|\leq\|x\|\ \forall x\in \R^m,$$  and  $a_i,\ a_{k+j}$ are non-zero constants. Without loss of generality, we will consider $a_1=1$.

Here, we consider the initial data on the Sobolev spaces, i.e., $$(U_0(x),V_0(x)):=((u_{0,1}, u_{0,2},...,u_{0,k}),(v_{0,1}, v_{0,2},...,v_{0,k}))\in H^{s}(\mathcal{G})=\bigoplus_{e \in E_-} H^s(-\infty, 0) \oplus\bigoplus_{e \in E_+} H^s(0, +\infty).$$

Our goal in studying the Cauchy problem \eqref{KDV}-\eqref{condfront} is to obtain results of local well-posedness, in Sobolev spaces with low regularity.

Initially, as mentioned previously, we denote the classical $L^2$ based Sobolev spaces on $\mathcal G$ by
\[ H^s(\mathcal{G}) = \bigoplus_{e \in E_-} H^s(-\infty, 0) \oplus \bigoplus_{e \in E_+} H^s(0, +\infty) .\]

We are interested in solving the problem in the set of regularity $s\in (-\frac12, \frac 32)\setminus \{\frac12\}$. Also, for $s\in (\frac12,\frac32)$, because of the existence of traces we consider the following compatibility condition
\begin{equation}\label{compatibility}
    u_{0,1} (0)= a_i u_{0,i} (0)= a_{k+j} v_{0,j} (0).
\end{equation}

\section{Statement of the main result} 

To introduce the main theorem, we need to establish some notation. In this sense, in this section, the index $i$ varies in integers between $1$ and $k$; $j$ integers between $1$ and $m$; and $l$ is $1$ or $2$. Moreover, in the statement of the theorem,
\begin{equation}\label{genmatrix}
\mathbf{M}(\boldsymbol\lambda, \boldsymbol\beta)=\left[\begin{array}{cc}
M_{1,1} & 0\\
M_{2,1} & M_{2,2}\\
M_{3,1} & M_{3,2}\\
M_{4,1} & M_{4,2}\\
\end{array}\right]
\end{equation}

\noindent is a $(2k+m) \times (2k+m)$ matrix formed by the following blocks, whose orders are shown subscript,\\[15pt]
\begin{equation}\label{m111}
M_{1,1}=\left[\begin{array}{ccccccccccccccccccccccccccc}
\rho_{1,1} & \rho_{1,2} & -a_{2}\rho_{2,1} & -a_{2}\rho_{2,2} & 0 & 0 & \cdots & 0 & 0\\
    \rho_{1,1} & \rho_{1,2} & 0 & 0 & -a_{3}\rho_{3,1}& -a_{3}\rho_{3,2} & \cdots & 0 & 0\\
    \vdots & \vdots & \vdots & \vdots & & & \ddots & \vdots & \vdots\\
    \rho_{1,1} & \rho_{1,2} & 0 & 0 & 0 & 0 & \cdots & -a_{k}\rho_{k,1} & -a_{k}\rho_{k,2}
\end{array}\right]_{(k-1) \times 2k},\\[5pt]
\end{equation}
\begin{equation}\label{m1212}
M_{2,1}=\left[\begin{array}{cccccccccccccccccccccccccccccc}
\rho_{1,1} & \rho_{1,2} & 0 & 0 & \hdots & 0\\
    \rho_{1,1} & \rho_{1,2} & 0 & 0 & \hdots & 0\\
    \rho_{1,1} & \rho_{1,2} & 0 & 0 & \hdots & 0\\
    \rho_{1,1} & \rho_{1,2} & 0 & 0 & \hdots & 0\\
\end{array}\right]_{m\times 2k},\\[5pt]
\end{equation}
\begin{equation}\label{m222}
M_{2,2}=\left[\begin{array}{ccccccccccccccccccccccccccccccccccc}
    -a_{k+1}d_1 & 0 & \cdots & 0\\
    0 & -a_{k+2}d_2 & \cdots & 0\\
    \vdots & \vdots & \ddots & \vdots\\
    0 & 0 & \cdots & -a_{k+m}d_m\\
\end{array}\right]_{m \times m},
\end{equation}

\begin{equation}\label{m3131}
M_{3,1}=\begin{bmatrix}
-\tilde{\rho}_{1,1} & -\tilde{\rho}_{1,2} & 0 & 0 & \cdots & 0 & 0 \\
0 & 0 & -\tilde{\rho}_{2,1} & -\tilde{\rho}_{2,2} & \cdots & 0 & 0 \\
\vdots & \vdots & \vdots & \vdots & \ddots & \vdots & \vdots \\
0 & 0 & 0 & 0 & \cdots & -\tilde{\rho}_{k,1} & \tilde{\rho}_{k,2} \\
\end{bmatrix}_{k \times 2k},
\end{equation}

\begin{equation}\label{m3232}
M_{3,2}=\begin{bmatrix}
b_{1,1} \tilde{d}_{1} & b_{1,2} \tilde{d}_{2} & \cdots & b_{1,m} \tilde{d}_{m} \\
b_{2,1} \tilde{d}_{1} & b_{2,2} \tilde{d}_{2} & \cdots & b_{2m} \tilde{d}_{m} \\
\vdots & \vdots & \ddots & \vdots \\
b_{k1} \tilde{d}_{1} & b_{k,2} \tilde{d}_{2} & \cdots & b_{k,m} \tilde{d}_{m}
\end{bmatrix}_{k \times m},
\end{equation}

\begin{equation}\label{m4141}
M_{4,1}=\begin{bmatrix}
a_{1}^{-1} \tilde{\tilde{\rho}}_{1,1} & a_{1}^{-1} \tilde{\tilde{\rho}}_{1,2} & a_{2}^{-1} \tilde{\tilde{\rho}}_{2,1} & a_{2}^{-1} \tilde{\tilde{\rho}}_{2,2} & \cdots & a_{k}^{-1} \tilde{\tilde{\rho}}_{k, 1} & a_{k}^{-1} \tilde{\tilde{\rho}}_{k, 2}
\end{bmatrix}_{1 \times 2k}
\end{equation}
and
\begin{equation}
M_{4,2}=\begin{bmatrix}\label{m4242}
-a_{k+1}^{-1}\tilde{\tilde{d}}_1 & -a_{k+2}^{-1}\tilde{\tilde{d}}_2 & \hdots & -a_{k+m}^{-1}\tilde{\tilde{d}}_m
\end{bmatrix}_{1 \times m}.
\end{equation}

All matrices, blocks, and terms $\rho_{i,l}; \rho_{k+j, l}; d_{j}; \tilde{\rho}_{il}; \tilde{\rho}_{k+j, l}; \tilde{d}_{j}; \tilde{\tilde{\rho}}_{il}; \tilde{\tilde{\rho}}_{k+j, l}; \tilde{\tilde{d}}_{j}$ above are functions of $(\lambda_{1,1},\hdots\lambda_{k,1},\lambda_{1,2},\hdots\lambda_{k,2},\beta_1,\hdots\beta_m)$; Table \ref{termosmatriz} shows the mathematical expressions of this dependence. The numbers $a_i,a_{k+j}$ and $b_{ij}$ come from the boundary conditions, with $b_{ij}$ being the generic term of the matrix $B$.

\begin{table}[H]
\begin{tabular}{cc}
 term & expression on parameters\\
 $\rho_{il}$ & $2\sin\left(\frac{\pi}{3}\lambda _{il}+\frac{\pi}{6}\right)$\\[4pt]

$d_j$ & $e^{i\pi\beta_j}$\\[4pt]

$\tilde \rho_{il}$ & $2\sin\left(\frac{\pi}{3}\lambda _{il}-\frac{\pi}{6}\right)$\\[4pt]

$\tilde{d}_j$ & $e^{i(\pi\beta_j-1)}$\\[4pt]

$\tilde{\tilde\rho}_{il}$ & $2\sin\left(\frac{\pi}{3}\lambda _{il}-\frac{\pi}{2}\right)$\\[4pt]

$\tilde{\tilde d}_j$ & $e^{i(\pi\beta_j-2)}$\\[4pt]
 
\end{tabular}
\caption{Expressions of matrix terms}
\label{termosmatriz}
\end{table}

Now, we can state the main result of this work, which gives a criterium to obtain a solution for the Cauchy problem \eqref{KDV} with boundary conditions $\eqref{condfront}$.
\begin{theorem}\label{grandeteorema}
    Assume $s\in(-\frac{1}{2},\frac{3}{2})\setminus\{\frac{1}{2}\}$ to be a fixed real number and $u_{0,i}\in H^s(-\infty,0)$ and $v_{0,j}\in H^s(0,+\infty)$. 
        Suppose that there exist $\lambda_{il}(s)$ and $\beta_j(s)$ in  that satisfy 
        \begin{equation}\label{hyp}
            \displaystyle max\{s-1,0\} < \lambda_{il}(s),\beta_j(s) < min\left\{s+\frac{1}{2},\frac{1}{2}\right\}
        \end{equation}
        and such that the matrix $\mathbf{M}(\boldsymbol\lambda, \boldsymbol\beta)$  is invertible. Then, there exist a time $T>0$, depending on the size of the initial data, and a solution $(u_1,\hdots,u_k,v_1,\hdots,v_m)$ in the space $C([0,T];H^s(\mathcal{G}))$  of the Cauchy problem $\eqref{KDV}$ with boundary conditions $\eqref{condfront}$ satisfying the compatibility conditions \eqref{compatibility} when $s>\frac{1}{2}$ (since the trace operator is well defined in $H^s(\R^+)$ for $s>\frac{1}{2}$). Furthermore, the data-to-solution map $(u_{0,1},\hdots,u_{0k},v_{0,1},\hdots,v_{0m})\mapsto (u_1,\hdots,u_k,v_1,\hdots,v_m)$ is locally Lipschitz continuous from $H^s(\mathcal{G})$ to $C([0, T];H^s(\mathcal{G}))$.
\end{theorem}

\begin{remark} [About the range of s] Firstly, we consider the regularity parameter s in the interval  $(-\frac12,\frac32)$ due to the limitations on the estimates for the Duhamel boundary operators classes in Lemma \ref{cof} (see also Remark \ref{remark}). Regarding the absence of the index $s=\frac12$, this is related to the understanding of compatibility conditions on the boundary. The approach used here is based on the method to solve nonlinear dispersive equations on the half-line by Colliander and Kenig \cite{CK}, and Holmer \cite{Holmer}. In this approach, the exact compatibility conditions on the boundary are crucial, as pointed out in the work \cite{Holmer}, at the level of regularity $s=\frac12$ the compatibility conditions are more complicated to formulate. 
Therefore, consistent with other works in the half-line context, we decided not to explore this regularity here.
\end{remark}

The approach used to prove the main result is based on the arguments introduced by Colliander and Kenig \cite{CK}  and improved by Holmer \cite{Holmer}, to solve nonlinear dispersive equations on half-lines based on the Riemann-Liouville fractional integral operator, and its adaptation to studying these equations in star graphs in \cite{fourth} and \cite{MC}. The main idea to prove Theorem \ref{grandeteorema} is the construction of an auxiliary forced Cauchy problem in all $\R$, analogous to \eqref{KDV}; more precisely:
\begin{equation}\label{KDVforcing}
\begin{cases}
\partial_t u_i+\partial^3_{x}u_i+u_i\partial_{x}u_i=\mathcal{T}_1^{i}(x)h_1^{i}(t)+\mathcal{T}_2^{i}(x)h_2^{i}(t),& (x, t)\in \R\times (0,T),\ i=1,2, ..., k;\\
\partial_t v_j+\partial^3_{x}v_j+v_j\partial_{x}v_j=\mathcal{T}_3^{j}(x)h_3^{j},& (x,t)\in \R\times (0,T),\ j=1,2, ..., m;\\
u_i(x,0)=\widetilde{u}_{0,i}(x),&  x \in  \R, i=1,2,..., k;\\
v_j(x,0)=\widetilde{v}_{0,j}(x),&  x \in  \R, j=1,2,..., m;\\
\end{cases}
\end{equation}
where $\mathcal{T}_1^{i}$ and $\mathcal{T}_2^{i}$ ($i=1,2,...,k$) are distributions supported on the positive half-line $\R^+$, $\mathcal{T}_3^{j}$ ($j=1,2,...,m$) are distributions supported on the negative half-line $\mathbb{R}^-$, $\widetilde{u}_{0,i}$ and $\widetilde{v}_{0,j}$ are nice extensions of $u_{0,i}$ and $v_{0,j}$ in $\mathbb{R}$. The unknown boundary  forcing functions $h_1^{i}$, $h_2^{i}$ and $h_3^{j}$ depend on the traces of unknown functions $u_i, v_j$ and their derivatives $\partial_x u_i, \partial_xv_j $, $\partial_x^2u_i$ and $\partial_x^2v_j$. These are selected to ensure that the vertex conditions are satisfied.

The solution of forced Cauchy problem \eqref{KDVforcing} satisfying the vertex conditions is constructed using the classical restricted norm method of Bourgain (see \cite{Bourgain1}) and the inversion of a Riemann-Liouville fractional integration operator.

Following \cite{MC} and \cite{Holmer} we consider the distributions $\mathcal{T}_1^{i}=\frac{x_{+}^{\lambda_{1,i}-1}}{\Gamma(\lambda_{1,i})}$, $\mathcal{T}_2^{i}=\frac{x_{+}^{\lambda_{2,i}-1}}{\Gamma(\lambda_{2,i})}$,  and $\mathcal{T}_3^{j}=\frac{x_{-}^{\lambda_{3,j}-1}}{\Gamma(\lambda_{3,j})}$  where 
\begin{equation}
\left\langle \frac{x_+^{\lambda-1}}{\Gamma(\lambda)},\phi\right\rangle=\int_0^{+\infty}\frac{x^{\lambda-1}}{\Gamma(\lambda)}\phi(x)dx,\; \text{for}\; \text{Re}\ \lambda>0.
\end{equation}
For other values of $\lambda$ we can define $\frac{x_+^{\lambda-1}}{\Gamma(\lambda)}=\frac{d^k}{dx}\frac{x^{\lambda+k-1}}{\Gamma(\lambda+k)}$, for any integer $k$ satisfying $k+\text{Re}\lambda>0$. Finally, we define $\frac{x_{-}^{\lambda-1}}{\Gamma(\lambda)}=e^{i\pi \lambda}\frac{(-x)_{+}^{\lambda-1}}{\Gamma(\lambda)}$.

The crucial point here are the appropriate choices of the parameters $\lambda_{1,i}$, $\lambda_{2,i}$, $\lambda_{3,j}$ and the search for unknown functions $h_1^{i}$, $h_2^{i}$, $h_3^{j}$, which will depend on the regularity index $s$.  

\section{Notations and function spaces} For $\phi=\phi(x)\in S(\mathbb{R})$,  $\displaystyle \hat{\phi}(\xi)=\int e^{-i\xi x}\phi(x)dx$ denotes the Fourier transform of $\phi$. For $u=u(x,t)\in S(\mathbb{R}^2)$, $$\hat{u}=\hat{u}(\xi,\tau)=\displaystyle \int e^{-i(\xi x+\tau t)}u(x,t)dxdt$$ denotes its space-time Fourier transform. We will denote the space and the time Fourier transform of $u(x,t)$, respectively by $\mathcal{F}_{x}u(\xi,t)$ and $\mathcal{F}_{t}u(x,\tau)$.  

For any real number $\xi$ we put $\langle \xi \rangle:=1+|\xi|$ and $f(\xi,\tau) \lesssim g(\xi,\tau)$ means that there is a constant $C$ such that $f(\xi,\tau) \leq Cg(\xi,\tau)$ for all $(\xi, \tau)\in \R^2$. 
The characteristic function of an arbitrary set $A$ is denoted by $\chi_{A}$. Throughout the paper, we fix a cutoff function $\psi \in C_0^{\infty}(\mathbb{R})$ such that
\begin{equation}\label{cutt}
\psi(t)=
\begin{cases}
1 & \text{if}\; |t|\le 1,\\
0 & \text{if}\; |t|\ge 2
\end{cases}
\end{equation}
and $\mathbb{R}^*=\mathbb{R}\setminus \{0\}.$ 

\subsection{Function Spaces}

The approach used in this work is based on the ideas of \cite{fourth} and \cite{MC} as well as some tools of Colliander and Kenig \cite{CK} and Holmer \cite{Holmer}. Let us recall some function spaces in the context of half-lines and the entire line.

Initially, for $s\geq 0$ we say that $\phi \in H^s(\mathbb{R}^+)$ if there exists $\tilde{\phi}\in H^s(\mathbb{R})$ such that
$\phi=\tilde{\phi}|_{\R+}$.  In this case we set $\|\phi\|_{H^s(\mathbb{R}^+)}:=\inf\limits_{\tilde{\phi}}\|\tilde{\phi}\|_{H^{s}(\mathbb{R})}$. For $s\geq 0$ define $$H_0^s(\mathbb{R}^+)=\Big\{\phi \in H^{s}(\mathbb{R}^+);\,\text{supp} (\phi) \subset[0,+\infty) \Big\}.$$ For $s<0$, define $H^s(\mathbb{R}^+)$ and $H_0^s(\mathbb{R}^+)$  as the dual space of $H_0^{-s}(\mathbb{R}^+)$ and  $H^{-s}(\mathbb{R}^+)$, respectively. 

Also, we define 
$$C_0^{\infty}(\mathbb{R}^+)=\Big\{\phi\in C^{\infty}(\mathbb{R});\, \text{supp}(\phi) \subset [0,+\infty)\Big\}$$
and $C_{0,c}^{\infty}(\mathbb{R}^+)$ as those members of $C_0^{\infty}(\mathbb{R}^+)$ with compact support. We recall that $C_{0,c}^{\infty}(\mathbb{R}^+)$ is dense in $H_0^s(\mathbb{R}^+)$ for all $s\in \mathbb{R}$. A definition for $H^s(\R^-)$ and $H_0^s(\R^-)$ can be given analogous to that for $H^s(\R^+)$ and $H_0^s(\R^+)$.

The following results summarize useful properties of the Sobolev spaces on the half-line. For the proofs we refer the reader to Section 2 of \cite{CK}.

\begin{lemma}\label{sobolevh0}
	For all $f\in H^s(\mathbb{R})$  with $-\frac{1}{2}<s<\frac{1}{2}$ we have
	\begin{equation*}	\|\chi_{(0,+\infty)}f\|_{H^s(\mathbb{R})}\lesssim  \|f\|_{H^s(\mathbb{R})}.
	\end{equation*}
\end{lemma}
\begin{lemma}\label{alta}
	If $\frac{1}{2}<s<\frac{3}{2}$ the following statements are valid:
	\begin{enumerate}
		\item [(a)] $H_0^s(\R^+)=\big\{f\in H^s(\R^+);f(0)=0\big\},$\medskip
		\item [(b)] If  $f\in H^s(\R^+)$ with $f(0)=0$, then $\|\chi_{(0,+\infty)}f\|_{H_0^s(\R^+)}\lesssim \|f\|_{H^s(\R^+)}$.
	\end{enumerate}
\end{lemma}

\begin{lemma}\label{cut}
	If $f\in  H_0^s(\mathbb{R}^+)$ with $s\in \R$, we then have
		\begin{equation*}
	\|\psi f\|_{H_0^s(\mathbb{R}^+)}\lesssim \|f\|_{H_0^s(\mathbb{R}^+)}.
	\end{equation*}
	
\end{lemma}
\begin{remark}
	In Lemmas \ref{sobolevh0}, \ref{alta} and \ref{cut} all constants $c$ depend only on $s$ and $\psi$.
\end{remark}

Now, we introduce adpated Bourgain spaces used by Holmer in the context of half-lines.

We denote by $X^{s,b}$ the so-called Bourgain spaces \cite{Bourgain1} associated with the linear KdV equation. More precisely, $X^{s,b}$  is the completion of $S(\mathbb{R}^2)$ with respect to the norm
\begin{equation*}\label{Bourgain-norm}
\|w\|_{X^{s,b}}=\|\langle\xi\rangle^s\langle\tau-\xi^3\rangle^b\hat{w}(\xi,\tau) \|_{L_{\tau}^2L^2_{\xi}}.
\end{equation*}
To obtain our results we also need to define the following auxiliary modified Bougain spaces of \cite{Holmer}. Let  $U^{s,b}$ and $V^{\alpha}$ be the completion of $S(\R^2)$ with respect to the norms:
\begin{align*}
&\|w\|_{U^{s,b}}=\left(\int\int \langle \tau\rangle^{2s/3} \langle \tau-\xi^3\rangle^{2b} |\widehat{w}(\xi,\tau)|^2d\xi d\tau\right)^{\frac{1}{2}}\label{bourgain-uxiliar-2}\\\intertext{and}
&\|w\|_{V^{\alpha}}=\left(\int\int  \langle \tau\rangle^{2\alpha} |\widehat{w}(\xi,\tau)|^2d\xi d\tau\right)^{\frac{1}{2}}.
\end{align*}

If $\Omega$ is a domain in $\mathbb{R}^2$, then define by $X^{s, b}(\Omega)$ and $Y^{s, b}(\Omega)$ restrictions of $X^{s, b}$ and $Y^{s, b}$ on $\Omega$, respectively, with natural restriction norms.

On the context of a star graph $\mathcal{G}$ we consider the Bourgain spaces $X^{s, b}(\mathcal{G} \times (0, T))$ by

\begin{equation}
\mathcal{X}^{s, b}(\mathcal{G} \times (0, T)) 
= \bigoplus_{e \in E_-} X^{s, b}((-\infty, 0) \times (0, T)) \oplus \bigoplus_{e \in E_+} X^{s, b}((0, \infty) \times (0, T)).
\end{equation}

In a similar way, we can define the spaces $\mathcal U^{s, b}(\mathcal{G} \times (0, T))$
 and $\mathcal V^{\alpha}(\mathcal{G} \times (0, T))$.

 For solving the Cauchy problem \eqref{KDV}-\eqref{condfront} we will consider the following functional space $\mathcal Z^{s, b,  \alpha}(T)$, $T > 0$, $s \ge 1$,

\begin{equation}\label{functional10}
    \mathcal Z^{s, b, \alpha}(T) = \left\{ w : \mathcal{G} \times [0, T] \to \mathbb{R} : 
    \begin{aligned}
        &w \in C([0, T]; H^s(\mathcal{G})), \\
        &w\in C(\mathcal{G}; H^{\frac{s+1}{3}}([0, T])) \cap \mathcal X^{s, b}(\mathcal G\times (0,T))\cap \mathcal V^{\alpha}(\mathcal G\times (0,T)), \\
        &\partial_x w \in C(\mathcal{G}; H^{\frac{s}{3}}([0, T])), \partial_x^2 w \in C(\mathcal{G}; H^{\frac{s-1}{3}}([0, T]))
    \end{aligned}
    \right\}.
\end{equation}
 
with the following norm
\begin{equation*}
\|(u_1,\hdots,u_k,v_1,\hdots,v_m)\|_{\mathcal Z^{(s,b,\alpha)}}:=\displaystyle\sum_{i=1}^k\|u_i\|_{\mathcal Z_i^{(s,b,\alpha)}}+\displaystyle\sum_{j=1}^m\|v_j\|_{\mathcal Z_{k+j}^{(s,b,\alpha)}},
\end{equation*}
where
\begin{equation}
\begin{split}
\|w\|_{\mathcal Z_i^{(s,b,\alpha)}}=&\|w\|_{C([0,T];H^s(\R^+))}+\|w\|_{C(\R^+;H^{\frac{s+1}{3}}([0,T]))}+\|w\|_{ X^{s,b}}+\|w\|_{V^{\alpha}}\\
&+\|w_x\|_{C(\R^+;H^{\frac{s}{3}}([0,T]))}+\|w_{xx}\|_{C(\R^+;H^{\frac{s-1}{3}}([0,T]))}.
\end{split}
\end{equation}

The next nonlinear estimates, in the context of the KdV equation for $b<\frac{1}{2}$, were derived by Holmer in Lemma 5.10 of \cite{Holmer}.
\begin{lemma}\label{bilinear1}
	\begin{itemize}
		\item [(a)]
		Given $s>-\frac{3}{4}$, there exists $b=b(s)<\frac12$ such that for all $\alpha>\frac{1}{2}$
		we have
		\begin{equation*}
		\big\|\partial_x (v_1v_2)\big\|_{X^{s,-b}}\lesssim \|v_1\|_{X^{s,b}\cap V^{\alpha}}\|v_2\|_{X^{s,b}\cap V^{\alpha}}.
		\end{equation*}
		\item[(b)] Given $-\frac{3}{4}<s<3$, there exists $b=b(s)<\frac12$ such that for all $\alpha>\frac{1}{2}$
		we have
		\begin{equation*}
		\big\|\partial_x (v_1v_2)\big\|_{Y^{s,-b}}\lesssim\|v_1\|_{X^{s,b}\cap V^{\alpha}}\|v_2\|_{X^{s,b}\cap V^{\alpha}}.
		\end{equation*}
	\end{itemize}
\end{lemma}

\subsection{Riemann-Liouville fractional integral}
For Re $\alpha>0$, the tempered distribution $\frac{t_+^{\alpha-1}}{\Gamma(\alpha)}$ is defined as a locally integrable function by 
\begin{equation*}
\left \langle \frac{t_+^{\alpha-1}}{\Gamma(\alpha)},\ f \right \rangle:=\frac{1}{\Gamma(\alpha)}\int_0^{+\infty} t^{\alpha-1}f(t)dt.
\end{equation*}

For Re $\alpha>0$, integration by parts implies that
\begin{equation*}
\frac{t_+^{\alpha-1}}{\Gamma(\alpha)}=\partial_t^k\left( \frac{t_+^{\alpha+k-1}}{\Gamma(\alpha+k)}\right)
\end{equation*}
for all $k\in\mathbb{N}$. This expression allows to extend the definition, in the sense of distributions,  of $\frac{t_+^{\alpha-1}}{\Gamma(\alpha)}$ to all $\alpha \in \mathbb{C}$.

If $f\in C_0^{\infty}(\mathbb{R}^+)$, we define
\begin{equation*}
\mathcal{I}_{\alpha}f=\frac{t_+^{\alpha-1}}{\Gamma(\alpha)}*f.
\end{equation*}
Thus, for Re $\alpha>0$,
\begin{equation*}
\mathcal{I}_{\alpha}f(t)=\frac{1}{\Gamma(\alpha)}\int_0^t(t-s)^{\alpha-1}f(s)ds
\end{equation*} 
and notice that 
$$\mathcal{I}_0f=f,\quad  \mathcal{I}_1f(t)=\int_0^tf(s)ds,\quad \mathcal{I}_{-1}f=f'\quad  \text{and}\quad  \mathcal{I}_{\alpha}\mathcal{I}_{\beta}=\mathcal{I}_{\alpha+\beta}.$$

The following result states important properties of the Riemann-Liouville fractional integral operators and their proof can be found in Lemma 5.4 of \cite{Holmer}.

\begin{lemma}\label{lio-lemaint}
	If $0\leq \alpha <\infty$,\, $s\in \mathbb{R}$ and $\varphi \in C_0^{\infty}(\mathbb{R})$, then we have
	\begin{align}
	&\|\mathcal{I}_{-\alpha}h\|_{H_0^s(\mathbb{R}^+)}\leq c \|h\|_{H_0^{s+\alpha}(\mathbb{R}^+)}\label{lio}\\
	\intertext{and}
	&\|\varphi \mathcal{I}_{\alpha}h\|_{H_0^s(\mathbb{R}^+)}\leq c_{\varphi} \|h\|_{H_0^{s-\alpha}(\mathbb{R}^+)}.\label{lemaint}
	\end{align}
\end{lemma}
\section{Estimates associated of some operators}\label{section3}

\subsection{Free popagator} The  unitary  group $e^{-t\partial_x^3}:S'(\mathbb{R})\rightarrow S'(\mathbb{R})$ associated to the linear KdV equation is defined by
\begin{equation*}
e^{-t\partial_x^3}\phi(x)=\Big(e^{it\xi^3}\widehat{\phi}(\xi)\Big)^{\lor{}}(x),
\end{equation*}
that satisfies
\begin{equation}\label{lineark}
\begin{cases}
(\partial_t+\partial_x^3)e^{-t\partial_x^3}\phi(x,t)=0& \text{for}\quad (x,t)\in \mathbb{R}\times\mathbb{R},\\
e^{-t\partial_x^3}(x,0)=\phi(x)&\text{for}\quad  x\in\mathbb{R}.
\end{cases} 
\end{equation}

We will establish a series of estimates for the solutions of linear IVP \eqref{lineark}. We refer the reader to \cite{Holmer} for a proof of these results.

\begin{lemma}\label{grupok}
	Let $s\in\mathbb{R}$ and $0< b<1$. If $\phi\in H^s(\mathbb{R})$, then we have \medskip
	\begin{enumerate}
		\item[(a)] (\textbf{space traces}) $\|e^{-t\partial_x^3}\phi(x)\|_{C\big(\mathbb{R}_t;\,H^s(\mathbb{R}_x)\big)}\lesssim \|\phi\|_{H^s(\mathbb{R})}$; \medskip
		\item[(b)] (\textbf{time traces})$\|\psi(t)\partial_x^j e^{-t\partial_x^3}\phi(x)\|_{C\big(\mathbb{R}_x;\,H^{(s+1-j)/3}(\mathbb{R}_t)\big)}\lesssim \|\phi\|_{H^s(\mathbb{R})}$ , for $j\in \mathbb{N}$; \medskip 
		\item [(c)] (\textbf{Bourgain spaces}) $\|\psi(t)e^{-t\partial_x^3}\phi(x)\|_{X^{s,b}\cap V^{\alpha}}\lesssim \|\phi\|_{H^s(\mathbb{R})}$ .
	\end{enumerate}
\end{lemma}

\subsection{The Duhamel boundary forcing operator associated to the linear KdV equation}\label{section4}
Now we give the properties of the Duhamel boundary forcing operator introduced in \cite{Holmer}, that is
\begin{equation}\label{lk}
\begin{split}
\mathcal{V}g(x,t)&=3\int_0^te^{-(t-t')\partial_x^3}\delta_0(x)\mathcal{I}_{-2/3}g(t')dt'\\
&=3\int_0^t A\left(\frac{x}{(t-t')^{1/3}}\right)\frac{\mathcal{I}_{-2/3}g(t')}{(t-t')^{1/3}}dt',
\end{split}
\end{equation}
defined for all $g\in C_0^{\infty}(\mathbb{R}^+)$ and $A$ denotes the Airy function 
$$A(x)=\frac{1}{2\pi}\int_{\xi}e^{ix\xi}e^{i\xi^3}d\xi.$$

From definition of $\mathcal{V}$  follows 
\begin{equation}\label{forcingk}
\begin{cases}
(\partial_t+\partial_x^3)\mathcal{V}g(x,t)=3\delta_0(x)\mathcal{I}_{-\frac{2}{3}}g(t)& \text{for}\quad (x,t)\in \mathbb{R}\times\mathbb{R},\\
\mathcal{V}g(x,0)=0& \text{for}\quad x\in\mathbb{R}.
\end{cases}
\end{equation}

The proof of the results exhibited in this section was shown in \cite{Holmer}.

\begin{lemma}\label{lemacrb1}
	Let $g\in C_0^{\infty}(\mathbb{R}^+)$ and consider a fixed time  $t\in[0,1]$. Then, 
	\begin{enumerate}
		\item[(a)]The functions $\mathcal{V}g(\cdot,t)$ and $\partial_x\mathcal{V}g(\cdot,t)$ are continuous in $x$ for all $x\in\mathbb{R}$;
		\item[(b)] The function $\partial_x^2\mathcal{V}g(x,t)$ is continuous in $x$ for all $x\neq 0$ and has a step discontinuity
		of size $3\mathcal{I}_{\frac{2}{3}}g(t)$ at $x=0$. Also, $\partial_x^2\mathcal{V}g(x,t)$ satisfies the spatial decay bounds
		$$
		|\partial_x^2\mathcal{V}g(x,t)|\leq c_k \|f\|_{H^{k+2}}\langle x\rangle^{-k} \quad \text{for all}\; k\geq 0.
		$$
	\end{enumerate}
\end{lemma}

Since $A(0)=\frac{1}{3\Gamma\big(\frac{2}{3}\big)}$ from  \eqref{lk} we have that $\mathcal{V}g(0,t)=g(t).$
It follows from \eqref{lk} (for details see Lemma 4.2 in \cite{CK}) the function 
\begin{equation}
v(x,t)=e^{-t\partial_x^3}\phi(x)+\mathcal{V}\big(g-e^{-t\partial_x^3}\phi\big|_{x=0}\big)(x,t),
\end{equation}
where $g\in C_0^{\infty}(\R^+)$ and $\phi \in S(\R)$, solves the linear problem
\begin{equation}\label{linearr}
\begin{cases}
(\partial_t+\partial_x^3)v(x,t)=0& \text{for}\quad (x,t)\in \R^*\times\mathbb{R},\\
v(x,0)=\phi(x)& \text{for}\quad x\in\mathbb{R},\\
v(0,t)=g(t)& \text{for}\quad t\in(0, +\infty),
\end{cases}
\end{equation}
in the sense of distributions.

Now, we consider the second boundary forcing operator, introduced by Holmer \cite{Holmer}, associated to the linear KdV equation:
\begin{equation}\label{secondf}
\mathcal{V}^{-1}g(x,t)=\partial_x\mathcal{V}\mathcal{I}_{\frac13}g(x,t)=3\int_0^tA'\left( \frac{x}{(t-t')^{1/3}}\right)\frac{\mathcal{I}_{-\frac13}g(t')}{(t-t')^{2/3}}dt'.
\end{equation}
From Lemma \ref{lemacrb1}, for all $g\in C_0^{\infty}(\mathbb{R}^+)$ the function $\mathcal{V}^{-1}g(x,t)$ is continuous in $x$ on $x\in\mathbb{R}$; moreover, using that $A'(0)=-\frac{1}{3\Gamma(\frac{1}{3})}$, we get the relation  $\mathcal{V}^{-1}g(0,t)=-g(t).$

Also, from \eqref{secondf} (see \cite{Holmer}) we get that
\begin{equation}
\begin{cases}
(\partial_t+\partial_x^3)\mathcal{V}^{-1}g(x,t)=3\delta_0'(x)\mathcal{I}_{-\frac{1}{3}}g(t)& \text{for}\quad (x,t)\in \mathbb{R}\times\mathbb{R},\\
\mathcal{V}^{-1}g(x,0)=0& \text{for}\quad x\in\mathbb{R},
\end{cases}
\end{equation}
in the sense of distributions.

Now, let us give regular functions $v_0(x)$, $g(t)$ and $h(t)$, then define 
\[ \Bigg[\begin{array}{c}
h_1 \vspace{0.3cm}\\
h_2 \end{array} \Bigg]:=\frac{1}{3}
\Bigg[\begin{array}{cc}
2  & -1 \vspace{0.3cm}\\
-1 & -1 \end{array} \Bigg]
\Bigg[\begin{array}{c}
g-e^{-\cdot\partial_x^3}v_0|_{x=0}\vspace{0.3cm}\\
\mathcal{I}_{\frac{1}{3}}\big(h-\partial_xe^{-\cdot\partial_x^3}v_0|_{x=0}\big)\end{array} \Bigg].\]
Holmer \cite{Holmer} proves that the  following function $v(x,t)=e^{-t\partial_x^3}v_0(x)+\mathcal{V}h_1(x,t)+\mathcal{V}^{-1}h_2(x,t)$ solves the linear problem 
\begin{equation}\label{linearl}
\begin{cases}
(\partial_t+\partial_x^3)v(x,t)=0& \text{for}\quad (x,t)\in \mathbb{R}^*\times\mathbb{R},\\
v(x,0)=v_0(x)& \text{for}\quad x\in\mathbb{R},\\
v(0,t)=g(t)& \text{for}\quad t\in \R ,\\
\lim\limits_{x \rightarrow 0^{-}}\partial_xv(x,t)=h(t)& \text{for}\quad t\in \R,
\end{cases}
\end{equation}
in the sense of distributions.

\subsection{The Duhamel Boundary Forcing Operator Classes associated to linear KdV equation}

To extend our results to a broader range of regularity (see Remark \ref{remark}), we utilize two classes of boundary forcing operators, generalizing the operators $\mathcal{V}$ and $\mathcal{V}^{ 
-1}$
  introduced by Holmer \cite{Holmer}.

Let $\lambda\in \mathbb{C}$ with 
$\text{Re}\,\lambda>-3$ and $g\in C_0^{\infty}(\mathbb{R}^+)$. Define the operators
\begin{equation*}
\mathcal{V}_{-}^{\lambda}g(x,t)=\left[\frac{x_+^{\lambda-1}}{\Gamma(\lambda)}*\mathcal{V}\big(\mathcal{I}_{-\frac{\lambda}{3}}g\big)(\cdot,t)   \right](x)
\end{equation*}
and
\begin{equation*}
\mathcal{V}_{+}^{\lambda}g(x,t)=\left[\frac{x_-^{\lambda-1}}{\Gamma(\lambda)}*\mathcal{V}\big(\mathcal{I}_{-\frac{\lambda}{3}}g\big)(\cdot,t)   \right](x),
\end{equation*}
with $\frac{x_{-}^{\lambda-1}}{\Gamma(\lambda)}=e^{i\pi \lambda}\frac{(-x)_{+}^{\lambda-1}}{\Gamma(\lambda)}$. Then, 
using \eqref{forcingk} we have that
\begin{equation*}
(\partial_t+\partial_x^3)\mathcal{V}_{-}^{\lambda}g(x,t)=3\frac{x_{+}^{\lambda-1}}{\Gamma(\lambda)}\mathcal{I}_{-\frac{2}{3}-\frac{\lambda}{3}}g(t)
\end{equation*}
and
\begin{equation*}
(\partial_t+\partial_x^3)\mathcal{V}_{+}^{\lambda}g(x,t)=3\frac{x_{-}^{\lambda-1}}{\Gamma(\lambda)}\mathcal{I}_{-\frac{2}{3}-\frac{\lambda}{3}}g(t).
\end{equation*}

The following lemmas state properties of the operators classes $\mathcal{V}_{\pm}^{\lambda}$. For the proofs, we refer the reader to Lemmas 3.1, 3.2 and 5.8 of \cite{Holmer}.
\begin{lemma}[\textbf{Spatial continuity and decay properties for} $\boldsymbol{\mathcal{V}_{\pm}^{\lambda}g(x,t)}$]\label{holmer1}
	Let $g\in C_0^{\infty}(\mathbb{R}^+)$. Then, we have
	\begin{equation*}
	\mathcal{V}_{\pm}^{\lambda-2}g=\partial_x^2\mathcal{V}^{\lambda}\mathcal{I}_{\frac{2}{3}}g,\quad \mathcal{V}_{\pm}^{\lambda-1}g=\partial_x\mathcal{V}^{\lambda}\mathcal{I}_{\frac{1}{3}}g\quad \text{and}\quad \mathcal{V}_{\pm}^{0}g=\mathcal{V}g.
	\end{equation*}
	Also, $\mathcal{V}_{\pm}^{-2}g(x,t)$ has a step discontinuity of size $3g(t)$ at $x=0$, otherwise for $x\neq 0$, $\mathcal{V}_{ \pm}^{-2}g(x,t)$ is continuous in $x$. For $\lambda>-2$, $\mathcal{V}_{\pm}^{\lambda}g(x,t)$ is continuous in $x$ for all $x\in\mathbb{R}$.
\end{lemma}

\begin{lemma}[\textbf{Values of} $\boldsymbol{\mathcal{V}_{\pm}^{\lambda}g(x,t)}$ \textbf{at} $\boldsymbol{x=0}$]\label{holmer2}
	For $\emph{Re}\,\lambda>-2$ and $g\in C_0^{\infty}(\R^+)$ we have 
	\begin{align*}
	&\mathcal{V}_{-}^{\lambda}g(0,t)=2\sin \left(\frac{\pi}{3}\lambda+\frac{\pi}{6}\right)g(t)\\ 
	\intertext{and}
	&\mathcal{V}_{+}^{\lambda}g(0,t)=e^{i\pi\lambda}g(t).
	\end{align*}
\end{lemma}
\begin{lemma}\label{cof}
	Let $s\in\mathbb{R}$.  The following estimates are ensured:\medskip
	\begin{enumerate}
		\item[(a)] (\textbf{space traces}) $\|\mathcal{V}_{\pm}^{\lambda}g(x,t)\|_{C\big(\mathbb{R}_t;\,H^s(\R_x)\big)}\lesssim \|g\|_{H_0^{(s+1)/3}(\mathbb{R}^+)}$ for all $s-\frac{5}{2}<\lambda<s+\frac{1}{2}$, $\lambda<\frac{1}{2}$  and $\emph{supp}(g)\subset[0,1]$.\medskip

		\item[(b)] (\textbf{time traces})
		$\|\psi(t)\partial_x^j\mathcal{V}_{\pm}^{\lambda}g(x,t)\|_{C\big(\mathbb{R}_x;\,H_0^{(s+1)/3}(\mathbb{R}_t^+)\big)}\lesssim c \|g\|_{H_0^{(s+1)/3}(\mathbb{R}^+)}$ for all $-2+j<\lambda<1+j$, for $j\in\{0,1,2\}$.\medskip

		\item[(c)] (\textbf{Bourgain spaces}) $
		\big\|\psi(t)\mathcal{V}_{\pm}^{\lambda}g(x,t)\big\|_{X^{s,b}\cap V^{\alpha}}\lesssim c \|g\|_{H_0^{(s+1)/3}(\mathbb{R}^+)}
		$ for all $s-1\leq \lambda<s+\frac{1}{2}$, $\lambda<\frac{1}{2}$, $\alpha\leq\frac{s-\lambda+2}{3}$ and  $0\leq b<\frac{1}{2}$.
	\end{enumerate}
\end{lemma} 

\begin{remark}\label{remark} The range of regularity $s$ such that the estimates of Lemma \ref{cof} hold varies depending on the index $\lambda$. An important point is that for each $s$ given, we will choose an appropriate $\lambda$ depending on $s$. For example, when $\lambda = 0$, the Bourgain space estimates are valid only for $-1/2 < s < 1$. This dependence on $\lambda$ is one of the reasons why we must consider the entire family of operators, rather than focusing solely on $\mathcal{V}$. We also note that by analyzing Lemma \ref{cof}, the values of $s$ for which there exists an appropriate $\lambda$ such that all claims of this lemma hold are contained in the set $(-\frac{1}{2}, \frac{3}{2})$.
\end{remark}

\subsection{{The Duhamel Inhomogeneous Solution Operator}}\label{sectionduhamel}

Now, we introduce the classical inhomogeneous solution operator $\mathcal{K}$ associated to the KdV equation given by
\begin{equation*}
\mathcal{K}w(x,t)=\int_0^te^{-(t-t')\partial_x^3}w(x,t')dt',
\end{equation*}
which satisfies
\begin{equation}\label{DK}
\begin{cases}
(\partial_t+\partial_x^3)\mathcal{K}w(x,t) =w(x,t)&\text{for}\quad  (x,t)\in\mathbb{R}\times\mathbb{R},\\
\mathcal{K}w(x,0) =0 & \text{for}\quad x\in\mathbb{R}.
\end{cases}
\end{equation}

The following lemma summarizes useful estimates for the Duhamel inhomogeneous solution operators  $\mathcal{K}$  that will be used later in the proof of the main results and its proof can be seen in \cite{Holmer}.

\begin{lemma}\label{duhamelk} For all $s\in\mathbb{R}$ we have the following estimates:\medskip
	\begin{enumerate}
		\item[(a)](\textbf{space traces}) Let $0<b<\frac{1}{2}$, then
		\begin{equation*}
		\|\psi(t)\mathcal{K}w(x,t)\|_{C\big(\mathbb{R}_t;\,H^s(\mathbb{R}_x)\big)}\lesssim\|w\|_{X^{s,-b}}.
		\end{equation*}\medskip
		
		\item[(b)](\textbf{time traces})
		Let $0<b<\frac{1}{2}$ and $j\in \{0,1,2\}$, then
		\begin{equation*}
		\left\|\psi(t)\partial_x^j\mathcal{K}w(x,t)\right\|_{C\big(\mathbb{R}_x;\,H^{(s+1)/3}(\mathbb{R}_t)\big)}\lesssim \begin{cases}
		\|w\|_{X^{s,-b}} & \text{if}\; -1+j\leq s \leq \frac{1}{2}+j,\\
		\|w\|_{X^{s,-b}}+\|w\|_{U^{s,-b}}& \text{for all}\;  s\in\mathbb{R}.
		\end{cases}
		\end{equation*}\medskip
				
		\item[(c)](\textbf{Bourgain spaces estimates})
		Let $0<b<\frac{1}{2}$ and $\alpha>1-b$, then
		\begin{equation*}
		\|\psi(t)\mathcal{K}w(x,t)\|_{X^{s,b}\cap V^{\alpha}}\lesssim \|w\|_{X^{s,-b}}.
		\end{equation*}
	\end{enumerate}
\end{lemma}

\section*{Notation}
In following sections, the index $i$ varies in integers between $1$ and $k$; $j$ integers between $1$ and $m$; and $l$ is $1$ or $2$.
\section{Start of the Proof of  Theorem \ref{grandeteorema}}\label{section6}
Here, we present the first part of the proof of the main result.

Firstly, we will obtain an integral equation that solves the Cauchy problem \eqref{KDV}.
We start rewriting the first vertex condition \eqref{condfront}  in terms of matrices:
\begin{equation}\label{m1}
\left[\begin{array}{cccccccc}
1&-a_2& 0&\cdots&0&0&\cdots&0\\
1& 0&-a_3&\cdots&0&0&\cdots&0\\
\vdots&\vdots&\vdots&\ddots&\vdots&\vdots&\ddots&\vdots\\
1&0&0&\cdots&-a_k&0&\cdots&0\\
1&0&0&\cdots&0&-a_{k+1}&\cdots&0\\
\vdots&\vdots&\vdots&\ddots&\vdots&\vdots&\ddots&\vdots\\
1&0&0&\cdots&0&0&\cdots&-a_{k+m}
\end{array}\right]_{(k+m-1) \times (k+m)}
\left[\begin{array}{c}
u_1(0,t)\\
\vdots\\
u_k(0,t)\\
v_1(0,t)\\
\vdots\\
v_m(0,t)
\end{array}\right]_{(k+m)\times 1}=0.
\end{equation}

Now, if we write $B$ as 
\begin{equation}\label{B}
\left[\begin{array}{cccc}
b_{1,1}& b_{1,2}& \cdots&b_{1,m}\\
b_{2,1}& b_{2,2}& \cdots&b_{2,m}\\
\vdots& \vdots& \ddots&\vdots\\
b_{k,1}& b_{k,2}& \cdots&b_{k,m}\\
\end{array}\right]_{k\times m},
\end{equation}

then the second condition on \eqref{condfront} becomes
\begin{equation}\label{m2}
\left[\begin{array}{ccccccc}
-1&0&\cdots&0&b_{1,1}&  \cdots&b_{1,m}\\
0&-1&\cdots&0&b_{2,1}& \cdots&b_{2,m}\\
\vdots&\vdots&\ddots&\vdots&  \vdots&\ddots&\vdots\\0&0&\cdots&-1&b_{k,1}& \cdots&b_{k,m}\\
\end{array}\right]_{k\times (k+m)}
\left[\begin{array}{c}
\partial_xu_1(0,t)\\
\vdots\\
\partial_xu_k(0,t)\\
\partial_xv_1(0,t)\\
\vdots\\
\partial_xv_m(0,t)\\
\end{array}\right]_{(k+m)\times 1}=0.
\end{equation}
Finally, the last condition is
\begin{equation}\label{m3}
\left[\begin{array}{cccccc}
a_1^{-1}& \cdots& a_k^{-1}&-a_{k+1}^{-1}& \cdots& -a_{k+m}^{-1}\\
\end{array}\right]_{1\times (k+m)}
\left[\begin{array}{c}
\partial^2_{x} u_1(0,t)\\
\vdots\\
\partial^2_{x} u_k(0,t)\\
\partial^2_{x} v_1(0,t)\\
\vdots\\
\partial^2_{x} v_m(0,t)\\
\end{array}\right]_{(k+m)\times 1}=0.
\end{equation}
Let $\widetilde{u}_{0,i}$ and $\widetilde{v}_{0,j}$ be nice extensions of ${u}_{0,i}$ and ${v}_{0,j}$, respectively satisfying
\begin{equation*}
\|\widetilde{u}_{0,i}\|_{H^s(\R)}\leq c\|{u}_{0,i}\|_{H^s(\R^{+})} \text{ and }\; \|\widetilde{v}_{0,j}\|_{H^s(\R)}\leq c\|{v}_{0,j}\|_{H^s(\R^{+})}.
\end{equation*}
Initially, based on the works \cite{MC} and \cite{fourth}, we look for solutions in the following form
\begin{align*}
&{u}_{i}(x,t)=\mathcal{V}_-^{\lambda_{i1}}\gamma_{i1}(x,t)+\mathcal{V}_-^{\lambda_{i2}}\gamma_{i2}(x,t)+F_i(x,t),\\
&{v}_{j}(x,t)=\mathcal{V}_+^{\beta_j}\theta_j(x,t)+G_j(x,t),
\end{align*}
where $\gamma_{i1}, \gamma_{i2}$ ($i=1,...,k$) and $\theta_j$ ($j=1,...,m$) are unknown functions that will be chosen later, and
\begin{align*}
&F_i(x,t)=e^{it\partial_x^3}\widetilde{u}_{0,i}+\mathcal{K}(u_i\partial_xu_i)(x,t),\\ 
&G_j(x,t)=e^{it\partial_x^3}\widetilde{v}_{0,j}+\mathcal{K}(v_j\partial_xv_j)(x,t),
\end{align*}
where $\mathcal{K}$ is the classical inhomogeneous solution operator associated with the KdV.

Assuming that $\gamma_{i1},\gamma_{i2}$ and $\theta_{j}$ belong to $C_{0}^{\infty}(\mathbb{R}^+)$  we can apply  Lemma \ref{holmer2} to obtain the following trace values
\begin{align}
&u_i(0,t)=2\sin\left(\frac{\pi}{3}\lambda_{i1}+\frac{\pi}{6}\right)\gamma_{i1}(t)+2\sin\left(\frac{\pi}{3}\lambda _{i2}+\frac{\pi}{6}\right)\gamma_{i2}(t)+F_i(0,t),\label{tec1}\\
&v_j(0,t)=e^{i\pi\beta_j}\theta_j(t)+G_j(0,t),\label{tec2}
\end{align}
where $\gamma_{i1}(t):=\gamma_{i1}(0,t)$, $\gamma_{i2}(t):=\gamma_{i2}(0,t)$ and $\theta_j(t):=\theta_j(0,t)$.

Now we calculate the traces of first derivative functions. By Lemmas \ref{holmer1} and \ref{holmer2} we see that
\begin{align}
&\partial_x u_i(0,t)=2\sin\left(\frac{\pi}{3}\lambda_{i1}-\frac{\pi}{6}\right)\mathcal{I}_{-\frac{1}{3}}\gamma_{i1}(t)+2\sin\left(\frac{\pi}{3}\lambda _{i2}-\frac{\pi}{6}\right)\mathcal{I}_{-\frac{1}{3}}\gamma_{i2}(t)+\partial_xF_i(0,t),\label{tec4}\\
&\partial_x v_j=e^{i\pi(\beta_j-1)}\mathcal{I}_{-\frac13}\theta_j(t)+\partial_xG_j(0,t),\label{tec5}
\end{align}
In the same way, we calculate the traces of second derivatives functions,
\begin{align}
&\partial^2_{x} u_i(0,t)=2\sin\left(\frac{\pi}{3}\lambda _{i1}-\frac{\pi}{2}\right)\mathcal{I}_{-\frac{2}{3}}\gamma_{i1}(t)+2\sin\left(\frac{\pi}{3}\lambda _{i2}-\frac{\pi}{2}\right)\gamma_{i2}(t)+\partial_x^2F_i(0,t),\label{tec7}\\
&\partial^2_{x} v_j(0,t)=e^{i\pi(\beta_j-2)}\mathcal{I}_{-\frac23}\theta_{j}(t)+\partial_x^2G_j(0,t),\label{tec8}
\end{align}
Note that, by Lemmas \ref{holmer1} and \ref{holmer2}, these calculations are valid for Re $\lambda>0$.

Using \eqref{tec1} and \eqref{tec2} we can rewrite \eqref{m1} as
\begin{equation}\label{m1A}
M_1
\left[\begin{array}{cc}
P& 0\\
0& D\\
\end{array}\right]_{(k+m)\times(2k+m)}\\
\left[\begin{array}{c}
\gamma_{1,1}\\
\gamma_{1,2}\\
\vdots\\
\gamma_{k1}\\
\gamma_{k2}\\
\theta_1\\
\vdots\\
\theta_m\\
\end{array}\right]_{(2k+m)\times 1}
=-M_1\left[\begin{array}{c}
F_1(0,t)\\
\vdots\\
F_k(0,t)\\
G_1(0,t)\\
\vdots\\
G_m(0,t)
\end{array}\right]_{(k+m)\times 1},
\end{equation}

where  $M_1$ is the $(k+m-1)\times (k+m)$ matrix on identity \eqref{m1},
\begin{equation}
P=
\left[\begin{array}{ccccccc}
\rho_{1,1}&\rho_{1,2}& 0&0& \hdots&0&0\\
0& 0&\rho_{2,1}&\rho_{2,2}& &\vdots&\vdots\\
\vdots&\vdots& & & \ddots&0&0\\
0&0&0&0& &\rho_{k,1}&\rho_{k,2}
\end{array}\right]_{(k\times 2k)},
\ D=
\left[\begin{array}{cccc}
d_{1}& 0 & \hdots& 0\\
0 & d_{2} &  & 0 \\
\vdots & & \ddots &  \\
0&0& &d_m
\end{array}\right]_{(m\times m)},
\end{equation}

$\rho_{il}=2\sin\left(\frac{\pi}{3}\lambda _{il}+\frac{\pi}{6}\right),\ l=1,2$ and $d_j=e^{i\pi\beta_j}$.\\

Applying  the operator  $\mathcal{I}_{\frac{1}{3}}$ on expressions \eqref{tec4} and \eqref{tec5} we can rewrite \eqref{m2} as 
\begin{equation}\label{m2A}
M_2
\left[\begin{array}{cc}
\tilde{P}& 0\\
0& \tilde D\\
\end{array}\right]_{(k+m)\times (2k+m)}\\
\left[\begin{array}{c}
\gamma_{1,1}\\
\gamma_{1,2}\\
\vdots\\
\gamma_{k1}\\
\gamma_{k2}\\
\theta_1\\
\vdots\\
\theta_m\\
\end{array}\right]_{(2k+m)\times 1}
=-M_2 \left[\begin{array}{c}
\partial_x \mathcal{I}_{\frac{1}{3}}F_1(0,t)\\
\vdots\\
\partial_x \mathcal{I}_{\frac{1}{3}}F_k(0,t)\\
\partial_x \mathcal{I}_{\frac{1}{3}}G_1(0,t)\\
\vdots\\
\partial_x \mathcal{I}_{\frac{1}{3}}G_m(0,t)
\end{array}\right]_{(k+m)\times 1},
\end{equation}
where  $M_2$ is the $k\times (k+m)$ matrix on identity \eqref{m2},
\begin{equation}
\tilde P=
\left[\begin{array}{ccccccc}
\tilde\rho_{1,1}&\tilde\rho_{1,2}& 0&0& \hdots&0&0\\
0& 0&\tilde\rho_{2,1}&\tilde\rho_{2,2}& &\vdots&\vdots\\
\vdots&\vdots& & & \ddots&0&0\\
0&0&0&0& &\tilde\rho_{k,1}&\tilde\rho_{k,2}
\end{array}\right]_{(k\times 2k)},
\ \tilde D=
\left[\begin{array}{cccc}
\tilde{d}_{1}& 0 & \hdots& 0\\
0 & \tilde{d}_{2} &  & 0 \\
\vdots & & \ddots &  \\
0&0& &\tilde{d}_m
\end{array}\right]_{(m\times m)},
\end{equation}

$\tilde \rho_{il}=2\sin\left(\frac{\pi}{3}\lambda _{il}-\frac{\pi}{6}\right),\ l=1,2$ and $\tilde d_j=e^{i(\pi\beta_j-1)}$.

Finally, we apply 
 the operator  $\mathcal{I}_{\frac{2}{3}}$ on expressions \eqref{tec7} and \eqref{tec8} to rewrite \eqref{m3} as 
\begin{equation}\label{m3A}
M_3
\left[\begin{array}{cc}
\tilde{\tilde P}& 0\\
0& \tilde{\tilde D}\\
\end{array}\right]_{(k+m)\times (2k+m)}\\
\left[\begin{array}{c}
\gamma_{1,1}\\
\gamma_{1,2}\\
\vdots\\
\gamma_{k1}\\
\gamma_{k2}\\
\theta_1\\
\vdots\\
\theta_m\\
\end{array}\right]_{2k+m}
=-M_3\left[\begin{array}{c}
\partial_x^2 \mathcal{I}_{\frac{2}{3}}F_1(0,t)\\
\vdots\\
\partial_x^2 \mathcal{I}_{\frac{2}{3}}F_k(0,t)\\
\partial_x^2 \mathcal{I}_{\frac{2}{3}}G_1(0,t)\\
\vdots\\
\partial_x^2 \mathcal{I}_{\frac{2}{3}}G_m(0,t)
\end{array}\right]_{(k+m)\times 1},
\end{equation}
where  $M_3$ is the $1\times (k+m)$ matrix on identity \eqref{m3},
\begin{equation}
\tilde {\tilde P}=
\left[\begin{array}{ccccccc}
\tilde{\tilde\rho}_{1,1}&\tilde{\tilde\rho}_{1,2}& 0&0& \hdots&0&0\\
0& 0&\tilde{\tilde\rho}_{2,1}&\tilde{\tilde\rho}_{2,2}& &\vdots&\vdots\\
\vdots&\vdots& & & \ddots&0&0\\
0&0&0&0& &\tilde{\tilde\rho}_{k,1}&\tilde{\tilde\rho}_{k,2}
\end{array}\right]_{(k\times 2k)},
\ \tilde {\tilde D}=
\left[\begin{array}{cccc}
\tilde{\tilde{{d}}}_{1}& 0 & \hdots& 0\\
0 & \tilde{\tilde{{d}}}_{2} &  & 0 \\
\vdots & & \ddots &  \\
0&0& &\tilde{\tilde{{d}}}_m
\end{array}\right]_{(m\times m)},
\end{equation}

$\tilde{\tilde\rho}_{il}=2\sin\left(\frac{\pi}{3}\lambda _{il}-\frac{\pi}{2}\right),\ l=1,2$ and $\tilde{\tilde d}_j=e^{i(\pi\beta_j-2)}$.\\

Now, multiplying the matrices in \eqref{m1A},  \eqref{m2A}, and  \eqref{m3A} we infer that

\begin{equation}\label{mm1final}
\begin{split}
&\left[\begin{array}{cc}
M_{1,1} & 0\\
M_{2,1} & M_{2,2}
\end{array}\right]_{(k+m-1)\times (2k+m)}
\left[\begin{array}{c}
\gamma_{1,1}\\
\gamma_{1,2}\\
\vdots\\
\gamma_{k1}\\
\gamma_{k2}\\
\theta_1\\
\vdots\\
\theta_m\\
\end{array}\right]_{(2k+m)\times 1}
=\left[\begin{array}{c}
-F_1(0,t)+a_2F_2(0,t)\\
-F_1(0,t)+a_3F_3(0,t)\\
\vdots\\
-F_1(0,t)+a_kF_k(0,t)\\
-F_1(0,t)+a_{k+1}G_1(0,t)\\
\vdots\\
-F_1(0,t)+a_{k+m}G_m(0,t)
\end{array}\right]_{(k+m-1)\times 1},
\end{split}
\end{equation}
in which $M_{1,1},\ M_{2,1}$ and $M_{2,2}$ are given by \eqref{m111}-\eqref{m222}.

Also,

\begin{equation}\label{mm2final}
\begin{split}
\left[\begin{array}{cc}
M_{3,1} & M_{3,2}\\
\end{array}\right]&_{k\times (2k+m)}
\left[\begin{array}{c}
\gamma_{1,1}\\
\gamma_{1,2}\\
\vdots\\
\gamma_{k1}\\
\gamma_{k2}\\
\theta_1\\
\vdots\\
\theta_m\\
\end{array}\right]_{(2k+m)\times 1}\\[15pt]
=&\left[\begin{array}{c}
\partial_x \mathcal{I}_{\frac{1}{3}}F_1(0,t)-b_{1,1}\partial_x \mathcal{I}_{\frac{1}{3}}G_1(0,t)-\hdots-b_{1m}\partial_x \mathcal{I}_{\frac{1}{3}}G_m(0,t)\\
\partial_x \mathcal{I}_{\frac{2}{3}}F_2(0,t)-b_{2,1}\partial_x \mathcal{I}_{\frac{1}{3}}G_1(0,t)-\hdots-b_{2m}\partial_x \mathcal{I}_{\frac{1}{3}}G_m(0,t)\\
\vdots\\
\partial_x \mathcal{I}_{\frac{1}{3}}F_k(0,t)-b_{k1}\partial_x \mathcal{I}_{\frac{1}{3}}G_1(0,t)-\hdots-b_{km}\partial_x \mathcal{I}_{\frac{1}{3}}G_m(0,t)
\end{array}\right]_{k \times 1},
\end{split}
\end{equation}
as well as
\begin{equation}\label{mm3final}
\begin{split}
\left[\begin{array}{cc}
M_{4,1} & M_{4,2}\\
\end{array}\right]_{1 \times (2k+m)}
&\left[\begin{array}{c}
\gamma_{1,1}\\
\gamma_{1,2}\\
\vdots\\
\gamma_{k1}\\
\gamma_{k2}\\
\theta_1\\
\vdots\\
\theta_m
\end{array}\right]_{(2k+m)\times 1}\\[15pt]
=&\left[\begin{array}{c}
-\displaystyle\sum_{i=1}^k a_i^{-1}\partial_x^2 \mathcal{I}_{\frac{2}{3}}F_i(0,t) + \sum_{j=1}^m a_{k+j}^{-1}\partial_x^2 \mathcal{I}_{\frac{2}{3}}G_j(0,t)
\end{array}\right]_{1\times 1},
\end{split}
\end{equation}
where $M_{3,1},\ M_{3,2}, M_{4,1}$ and $M_{4,2}$ are the matrices in $\eqref{m3131}$-$\eqref{m4242}$.

From \eqref{mm1final}, \eqref{mm2final} and \eqref{mm3final} we need to get functions $\gamma_{il}$ and $\theta_j$, and parameters $\lambda_{il}$ and $\beta_j$ satisfying
\begin{equation}\label{M}
\begin{split}
\left[\begin{array}{cc}
M_{1,1} & 0\\
M_{2,1} & M_{2,2}\\
M_{3,1} & M_{3,2}\\
M_{4,1} & M_{4,2}\\
\end{array}\right]&_{(2k+m)\times (2k+m)}\left[\begin{array}{r}
\gamma_{1,1}\\
\gamma_{1,2}\\
\vdots\\
\gamma_{k1}\\
\gamma_{k2}\\
\theta_1\\
\vdots\\
\theta_m\\
\end{array}\right]_{(2k+m)\times 1} \\[15pt]=&\left[\begin{array}{c}
-F_1(0,t)+a_2F_2(0,t)\\
-F_1(0,t)+a_3F_3(0,t)\\
\vdots\\
-F_1(0,t)+a_kF_k(0,t)\\
-F_1(0,t)+a_{k+1}G_1(0,t)\\
\vdots\\
-F_1(0,t)+a_{k+m}G_m(0,t)\\
\partial_x \mathcal{I}_{\frac{1}{3}}F_1(0,t)-b_{1,1}\partial_x \mathcal{I}_{\frac{1}{3}}G_1(0,t)-\hdots-b_{1m}\partial_x \mathcal{I}_{\frac{1}{3}}G_m(0,t)\\
\partial_x \mathcal{I}_{\frac{2}{3}}F_2(0,t)-b_{2,1}\partial_x \mathcal{I}_{\frac{1}{3}}G_1(0,t)-\hdots-b_{2m}\partial_x \mathcal{I}_{\frac{1}{3}}G_m(0,t)\\
\vdots\\
\partial_x \mathcal{I}_{\frac{1}{3}}F_k(0,t)-b_{k1}\partial_x \mathcal{I}_{\frac{1}{3}}G_1(0,t)-\hdots-b_{km}\partial_x \mathcal{I}_{\frac{1}{3}}G_m(0,t)\\
-\displaystyle\sum_{i=1}^k a_i^{-1}\partial_x^2 \mathcal{I}_{\frac{2}{3}}F_i(0,t) + \sum_{j=1}^m a_{k+j}^{-1}\partial_x^2 \mathcal{I}_{\frac{2}{3}}G_j(0,t)
\end{array}\right]_{(2k+m)\times 1},\\
&
\end{split}
\end{equation}
since that $M_{1,1},\ M_{2,1},\ M_{2,2},\ M_{3,1},\ M_{3,2},\ M_{4,1} $ and $M_{4,2}$ depend on $\rho_{il},\ d_j,\ \tilde{\rho}_{il},\ \tilde{d}_j,\ \tilde{\tilde{\rho}}_{il}\text{ and } \tilde{\tilde{d}}_j$, which are functions on $\lambda_{il}$'s and $\beta_j$'s.

We denote a simplified notation of \eqref{M} as
\begin{equation}
\mathbf{M}(\boldsymbol\lambda, \boldsymbol\beta)\left[\begin{array}{c}
\boldsymbol\gamma\\
\boldsymbol\theta
\end{array}\right]
=\mathbf{F},
\end{equation}
where $\mathbf{M}(\boldsymbol\lambda, \boldsymbol\beta)$ is the first matrix that appears in \eqref{M}, $\displaystyle\left[\begin{array}{c}
\boldsymbol\gamma\\
\boldsymbol\theta
\end{array}\right]$ is the  matrix column given by vector $(\gamma_{1,1}, \gamma_{1,2}, \hdots, \gamma_{k1}, \gamma_{k2}, \theta_1, \hdots, \theta_m)$ and $\mathbf{F}$ is the last matrix in $\eqref{M}$. By using the hypothesis of Theorem \ref{grandeteorema} we fix parameters $\lambda_{il}$ and $\beta_j$, such that
 \begin{equation}
\max\{s-1,0\}<\lambda_{il}(s),\beta_j(s)<\min\left\{s+\frac12,\frac12\right\}.
 \end{equation}
and the matrix $\mathbf{M}(\boldsymbol\lambda, \boldsymbol\beta)$ is invertible.

\section{Truncated integral operator}

In this step of the proof, we will define the truncated integral operator and the appropriate function space.

Given $s$ as in the hypothesis of Theorem \ref{grandeteorema} we fix the parameters $\lambda_{il}$ and $\beta_j$ as well as the functions $\gamma_{il}$ and $\theta_{j}$ chosen as in section \ref{section6}. Let $b=b(s)<\frac12$ and $\alpha(b)>1/2$ be such that the estimates given in Lemma \ref{bilinear1} are valid.

Define the operator
\begin{equation}
\Lambda=(\Lambda^-_1,\hdots,\Lambda^-_k,\Lambda^+_1,\hdots,\Lambda^+_m)
\end{equation}
where
\begin{align*}
&\Lambda^-_i u_i(x,t)=\psi(t)\mathcal{V}_-^{\lambda_{i1}}\gamma_{i1}(x,t)+\psi(t)\mathcal{V}_-^{\lambda_{i2}}\gamma_{i2}(x,t)+F_i(x,t),\\
&\Lambda^+_jv_j(x,t)=\psi(t)\mathcal{V}_+^{\beta_j}\theta_j(x,t)+G_j(x,t),
\end{align*}
where
\begin{align*}
&F_i(x,t)=\psi(t)(e^{it\partial_x^3}\widetilde{u}_{0,i}+\mathcal{K}(u_i(u_i)_x)(x,t)),\\ &G_j(x,t)=\psi(t)(e^{it\partial_x^3}\widetilde{v}_{0,j}+\mathcal{K}(v_j(v_j)_x)(x,t)).
\end{align*}

Here, we recall the cutoff function $\psi$ defined in \eqref{cutt}, which is essential to use the Fourier restriction method introduced by Bourgain in \cite{Bourgain1}.

We consider $\Lambda$ on the Banach space $\displaystyle Z_{(s,b,\alpha)}=\displaystyle\prod_{i=1}^k Z_{(s,b,\alpha)}^{i}\cdot\prod_{j=1}^m Z^{k+j}_{(s,b,\alpha)}$, where
\begin{equation*}
\begin{split}
Z^n_{(s,b,\alpha)}=\{z\in C(\R_t;&H^s(\R_x)) \cap C(\R_x;H^{\frac{s+1}{3}}(\R_t))\cap X^{s,b}\cap V^{\alpha};\\
& z_x\in C(\R_x;H^{\frac{s}{3}}(\R_t)),z_{xx}\in C(\R_x;H^{\frac{s-1}{3}}(\R_t)) \}\; (n=1,\hdots,k+m),
\end{split}
\end{equation*}
with norm
\begin{equation*}
\|(u_1,\hdots,u_k,v_1,\hdots,v_m)\|_{Z_{(s,b,\alpha)}}=\displaystyle\sum_{i=1}^k\|u_i\|_{Z^i_{(s,b,\alpha)}}+\displaystyle\sum_{j=1}^m\|v_j\|_{Z^{k+j}_{(s,b,\alpha)}},
\end{equation*}
and
\begin{equation}
\begin{split}
\|z\|_{Z^{n}_{(s,b,\alpha)}}=&\|z\|_{C(\R_t;H^s(\R_x))}+\|z\|_{C(\R_x;H^{\frac{s+1}{3}}(\R_t))}+\|z\|_{ X^{s,b}}+\|z\|_{V^{\alpha}}\\
&+\|z_x\|_{C(\R_x;H^{\frac{s}{3}}(\R_t))}+\|z_{xx}\|_{C(\R_x;H^{\frac{s-1}{3}}(\R_t))}.
\end{split}
\end{equation}

\begin{remark}
 Note that the spaces $Z_{(s,b,\alpha)}$ defined previously consists of the functions with domain $\prod_{n=1}^{k+m}\mathbb R\times (0,T)$, unlike  the space $\mathcal Z^{s,b,\alpha}$ defined in \eqref{functional10}, where the functions have domain on $\mathcal G \times (0,T)$. 
\end{remark}

{Now, we will prove that the functions} ${\mathcal{V}_{-}^{\lambda_{il}}\gamma_{il}(x,t)}$ {and} $
{\mathcal{V}_{+}^{\beta_j}\theta_j(x,t)}$ {are well defined}.

By Lemma \ref{cof} it suffices to show that the functions $\gamma_{il}$ and $\theta_j$ are in the closure of the spaces $C_0^{\infty}(\R^+)$. By using expression \eqref{M} we see that the functions $\gamma_{il}$ and $\theta_j$ are linear combinations of the functions $F_1(0,t)-a_iF_i(0,t),\ F_1(0,t)-a_{k+j}G_{j}(0,t),\ \partial_x \mathcal{I}_{\frac{1}{3}}F_i(0,t)-b_{i1}\partial_x \mathcal{I}_{\frac{1}{3}}G_1(0,t)-\hdots-b_{im}\partial_x \mathcal{I}_{\frac{1}{3}}G_m(0,t)$ and $
-\sum_{i=1}^k a_i^{-1}\partial_x^2 \mathcal{I}_{\frac{2}{3}}F_i(0,t) + \sum_{j=1}^m a_{k+j}^{-1}\partial_x^2 \mathcal{I}_{\frac{2}{3}}G_m(0,t)$. Thus, we need to show that the functions $F_i(0,t),\ G_j(0,t),\  \partial_x\mathcal{I}_{\frac{1}{3}}F_i(0,t),\ \partial_x\mathcal{I}_{\frac{1}{3}}G_j(0,t)$, $\partial_x^2\mathcal{I}_{\frac{2}{3}}F_i(0,t)$ and $\partial_x^2\mathcal{I}_{\frac{2}{3}}G_j(0,t)$ are in appropriate spaces. By using Lemmas  \ref{grupok}, \ref{cof}, \ref{duhamelk} and \ref{bilinear1} we obtain
\begin{equation}
\|F_i(0,t)\|_{H^{\frac{s+1}{3}}(\R^+)}\leq c( \|u_{0,i}\|_{H^s(\R^+)}+\|u_{i}\|_{X^{s,b}}^2+\|u_{i}\|_{V^{\alpha}}^2),
\end{equation}
\begin{equation}
\|G_j(0,t)\|_{H^{\frac{s+1}{3}}(\R^+)}\leq c( \|v_{0,j}\|_{H^s(\R^+)}+\|v_{j}\|_{X^{s,b}}^2+\|v_{j}\|_{V^{\alpha}}^2).
\end{equation}
If $-\frac{1}{2}<s<\frac{1}{2}$ we have that $\frac16 <\frac{s+1}{3}<\frac12$. Thus Lemma \ref{sobolevh0} implies that $H^{\frac{s+1}{3}}(\R^+)=H_0^{\frac{s+1}{3}}(\R^+)$. It follows that $F_i(0,t),G_j(0,t)\in H_0^{\frac{s+1}{3}}(\R^+)$ for $-\frac{1}{2}<s<\frac{1}{2}.$

If $\frac{1}{2}<s<\frac{3}{2}$, then $\frac12<\frac{s+1}{3}<\frac56$. Using the compatibility condition \eqref{compatibility} in Theorem \ref{grandeteorema}, we have that 
\begin{equation*}
F_1(0,0)-a_i F_i(0,0)=u_1(0,0)-a_iu_i(0,0)=u_{0,1}(0)-a_iu_{0,i}(0)=0,
\end{equation*}
\begin{equation*}
F_1(0,0)-a_{k+j} G_j(0,0)=u_1(0,0)-a_{k+j}v_j(0,0)=u_{0,1}(0)-a_{k+j}v_{0,j}(0)=0.
\end{equation*}
Then Lemma \ref{alta} implies 
\begin{equation}\label{trace1}
\begin{split}
F_1(0,t)-a_iF_i(0,t)\in H_0^{\frac{s+1}{3}}(\R^+),\\
F_1(0,t)-a_{k+j}G_j(0,t)\in H_0^{\frac{s+1}{3}}(\R^+)
\end{split}
\end{equation}

Now, using Lemmas  \ref{grupok}, \ref{cof}, \ref{duhamelk} and \ref{bilinear1} we see that
\begin{equation*}
\|\partial_x F_i(0,t)\|_{H^{\frac{s}{3}}(\R^+)}\leq c( \|u_{0,i}\|_{H^s(\R^+)}+\|u_i\|_{X^{s,b}}^2+\|u_i\|_{V^{\alpha}}^2), 
\end{equation*}
\begin{equation*}
\|\partial_xG_j(0,t)\|_{H^{\frac{s}{3}}(\R^+)}\leq c( \|v_{0,j}\|_{H^s(\R^+)}+\|v_j\|_{X^{s,b}}^2+\|v_j\|_{V^{\alpha}}^2). 
\end{equation*}
Since $-\frac12<s<\frac32$ we have $-\frac16<\frac{s}{3}<\frac12$, then Lemma \ref{sobolevh0} implies that the functions $\partial_xF_i(0,t),\partial_xG_j(0,t)\in H_0^{\frac{s}{3}}(\R^+)$. Then using Lemma \ref{lio} we have that
\begin{equation*}
\|\partial_x \mathcal{I}_{\frac13}F_i(0,t)\|_{H_0^{\frac{s+1}{3}}(\R^+)}\leq c( \|u_{0,i}\|_{H^s(\R^+)}+\|u_i\|_{X^{s,b}}^2+\|u_i\|_{V^{\alpha}}^2), 
\end{equation*}
\begin{equation*}
\|\partial_x\mathcal{I}_{\frac13}G_j(0,t)\|_{H_0^{\frac{s+1}{3}}(\R^+)}\leq c( \|v_{0,j}\|_{H^s(\R^+)}+\|v_{j}\|_{X^{s,b}}^2+\|v_{j}\|_{V^{\alpha}}^2). 
\end{equation*}

Thus, we have 
\begin{equation}\label{trace2}
\partial_x\mathcal{I}_{\frac{1}{3}}F_i(0,t)-b_{i1}\partial_x\mathcal{I}_{\frac{1}{3}}G_1(0,t)-\hdots-b_{im}\partial_x\mathcal{I}_{\frac{1}{3}}G_m(0,t)\in H_0^{\frac{s+1}{3}}(\R^+).
\end{equation}

In the same way we can obtain
\begin{equation*}
\|\partial_x^2\mathcal{I}_{\frac{2}{3}}F_i(0,t)\|_{H_0^{\frac{s+1}{3}}(\R^+)}\leq c( \|u_{0,i}\|_{H^s(\R^+)}+\|u_i\|_{X^{s,b}}^2+\|u_i\|_{V^{\alpha}}^2), 
\end{equation*}
\begin{equation*}
\|\partial_x^2\mathcal{I}_{\frac{2}{3}}G_j(0,t)\|_{H_0^{\frac{s+1}{3}}(\R^+)}\leq c( \|v_{0,j}\|_{H^s(\R^+)}+\|v_j\|_{X^{s,b}}^2+\|v_j\|_{V^{\alpha}}^2). 
\end{equation*}

It follows that
\begin{equation}\label{trace3}
-\displaystyle\sum_{i=1}^k a_i^{-1}\partial_x^2 \mathcal{I}_{\frac{2}{3}}F_i(0,t) + \sum_{j=1}^m a_{k+j}^{-1}\partial_x^2 \mathcal{I}_{\frac{2}{3}}G_j(0,t)\in H_0^{\frac{s+1}{3}}(\R^+).
\end{equation}

Thus, \eqref{trace1}, \eqref{trace2} and \eqref{trace3} imply that the functions ${\mathcal{V}_{-}^{\lambda_{il}}\gamma_{il}(x,t)}$ {and} ${\mathcal{V}_{+}^{\beta_j}\theta_j(x,t)}$ are well defined.

\section{End of the proof: Obtaining a fixed point of $\Lambda$}
Using Lemmas  \ref{lio}, \ref{grupok}, \ref{cof}, \ref{duhamelk} and \eqref{bilinear1} we obtain
\begin{equation}
\begin{split}
\|\Lambda (\tilde{u}_{1},\hdots,\tilde{u}_{k},\tilde{v}_{1},\hdots,\tilde{v}_{m})-\Lambda (\tilde{\tilde{u}}_{1},\hdots,\tilde{\tilde{u}}_{k},\tilde{\tilde{v}}_{1},\hdots,\tilde{\tilde{v}}_{m})\|_{Z_{(s,b,\alpha)}}\\
\leq c(\| (\tilde{u}_{1},\hdots,\tilde{u}_{k},\tilde{v}_{1},\hdots,\tilde{v}_{m})\|_{Z_{(s,b,\alpha)}}+\| (\tilde{\tilde{u}}_{1},\hdots,\tilde{\tilde{u}}_{k},\tilde{\tilde{v}}_{1},\hdots,\tilde{\tilde{v}}_{m})\|_{Z_{(s,b,\alpha)}})\\ \cdot\| (\tilde{u}_{1},\hdots,\tilde{u}_{k},\tilde{v}_{1},\hdots,\tilde{v}_{m})- (\tilde{\tilde{u}}_{1},\hdots,\tilde{\tilde{u}}_{k},\tilde{\tilde{v}}_{1},\hdots,\tilde{\tilde{v}}_{m})\|_{Z_{(s,b,\alpha)}}
\end{split}
\end{equation}
and 
\begin{equation}
\begin{split}
\|\Lambda ({u}_{1},\hdots,{u}_{k},{v}_{1},\hdots,{v}_{m})\|_{Z_{(s,b,\alpha)}}&\leq c \left( \sum_{i=1}^k\|u_{0,i}\|_{H^s(\R^+)}+\sum_{j=1}^m\| v_{0,j}\|_{H^s(\R^+)}\right.\\
&\left.+\sum_{i=1}^k(\|u_i\|_{X^{s,b}}^2+\|u_i\|_{V^{\alpha}}^2)+\sum_{j=1}^m(\|v_j\|_{X^{s,b}}^2+\|v_j\|_{V^{\alpha}}^2)\right) .
\end{split}
\end{equation}

Taking $\sum_{i=1}^{k}\|u_{0,i}\|_{H^s(\R^+)}+\sum_{j=1}^{m}\|v_{0,j}\|_{H^s(\R^+)}<\delta$ for $\delta>0$ suitably small, by the Banach fixed-point theorem, we obtain a fixed point $\Lambda (\overline{u}_{1},\hdots,\overline{u}_{k},\overline{v}_{1},\hdots,\overline{v}_{m})=(\overline{u}_{1},\hdots,\overline{u}_{k},\overline{v}_{1},\hdots,\overline{v}_{m})$ in a ball (since $\Lambda$ is a contraction in that space) $$B=\{({u}_{1},\hdots,{u}_{k},{v}_{1},\hdots,{v}_{m})\in Z_{(s,b,\alpha)}; \|({u}_{1},\hdots,{u}_{k},{v}_{1},\hdots,{v}_{m})\|_{Z_{(s,b,\alpha)}}\leq 2c\delta\}.$$ It follows that the restriction
\begin{equation}
({u}_{1},\hdots,{u}_{k},{v}_{1},\hdots,{v}_{m})=(\overline{u_1}\big|_{\R^-\times (0,1)},\hdots,\overline{u_k}\big|_{\R^-\times (0,1)},\overline{v_1}\big|_{\R^+\times (0,1)},\hdots, \overline{v_m}\big|_{\R^+\times (0,1)})
\end{equation}
solves the Cauchy problem \eqref{KDV} with boundary conditions \eqref{condfront}, on the space $\mathcal Z_{s,b,\alpha}$, in the sense of distributions.

Finally, the existence of solutions for any data in $H^s(\mathcal{G})$ follows by the standard scaling argument. Suppose that we are given data $\widetilde{u}_{0,i}$ and $\widetilde{v}_{0,j}$ with arbitrary size for the Cauchy problem \eqref{KDV} with boundary conditions \eqref{condfront}. For $\lambda<<1$ (to be selected later), define $u_{0,i}(x)=\lambda^{2}\widetilde{u}_{0,i}(\lambda x)$ and $v_{0,i}(x)=\lambda^{2}\widetilde{v}_{0,i}(\lambda x)$. Taking $\lambda$ sufficiently small so that
\begin{equation}
\sum_{i=1}^k\|u_{0,i}\|_{H^s(\R^+)}+\sum_{j=1}^m\|v_{0,j}\|_{H^s(\R^+)}\leq \lambda^{\frac32}(1+\lambda^{s})\left(\sum_{i=1}^k\|\widetilde{u}_{0,i}\|_{H^s(\R^+)}+\sum_{j=1}^m\|\widetilde{v}_{0,j}\|_{H^s(\R^+)}\right)<\delta.
\end{equation}
Then using the previous argument, there exists a solution $(u_1,\hdots,u_k,v_1,\hdots,v_m)$ for the problem \eqref{KDV}, with boundary conditions \eqref{condfront}, on $0\leq t\leq 1$. Then $(\widetilde u_1,\hdots,\widetilde u_k,\widetilde v_1,\hdots,\widetilde v_m)$ with $\widetilde{u}_i(x,t)=\lambda^{-2}u(\lambda^{-1}x,\lambda^{-3}t)$ and $\widetilde{v}_i(x,t)=\lambda^{-2}v(\lambda^{-1}x,\lambda^{-3}t)$ solves the Cauchy problem for the initial data $\widetilde{u}_{0,i}$ and $\widetilde{v}_{0,j}$ on time interval $0\leq t\leq \lambda^3$.
\section{Applications of Theorem \ref{grandeteorema}}
In this section, we give applications of the general criteria described in Theorem \ref{grandeteorema} in some situations, including balanced and unbalanced star graphs.

\begin{example}[$\mathcal{Y}$-junction] In this first example, we consider the case where the parameters which appear in boundary conditions are given by $a_1=a_{2}=a_3=1$ and $B=\left[\frac{\sqrt{2}}{2}\quad \frac{\sqrt{2}}{2}\right]$ and the star graph with one edge $(-\infty, 0)$ and two edges $(0, \infty)$ $(k=1,\ m=2)$ is called $\mathcal{Y}$-junction. Using Theorem $\ref{grandeteorema}$, for each {$s\in \left(-0.06, 1.22\right)\setminus\left\{\frac{1}{2}\right\}$}, $\lambda_{il}$ and $\beta_j$ can be found  in order to make the matrix $\mathbf{M}(\boldsymbol\lambda, \boldsymbol\beta)$ invertible . There exists a time $T>0$ and a solution of the Cauchy problem \eqref{KDV}-\eqref{condfront}  satisfying the compatibility conditions \eqref{compatibility} when $s>\frac{1}{2}$. 
    
\end{example}

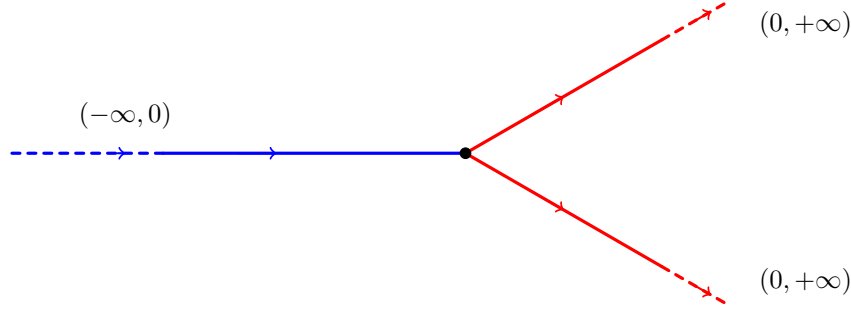
\begin{figure}[h]
    \centering
    \begin{tikzpicture}[line cap=round]
        \draw[dashed, blue, very thick] (-6,0) -- (-4,0);
        \draw[thick, blue, very thick] (-4,0) -- (0,0);
        \draw[blue, ->, thick] (-4.6,0) -- (-4.5,0);
        \draw[blue, ->, thick] (-3,0) -- (-2.5,0);
        \draw[thick, red, very thick] (0,0) -- (30:3);
        \draw[dashed, red, very thick] (30:3) -- (30:4);
        \draw[red, ->, thick] (30:1) -- (30:1.5);
        \draw[red, ->, thick] (30:3.5) -- (30:3.75);
        \draw[thick, red, very thick] (0,0) -- (-30:3);
        \draw[dashed, red, very thick] (-30:3) -- (-30:4);
        \draw[red, ->, thick] (-30:1) -- (-30:1.5);
        \draw[red, ->, thick] (-30:3.5) -- (-30:3.75);
        \node at (-4.5,0.5) {$( -\infty, 0)$};
        \node at (4.5,1.7) {$(0, +\infty)$};
        \node at (4.5,-1.7) {$(0, +\infty)$};
        \filldraw (0,0) circle (2pt);

    \end{tikzpicture}
    \caption{$\mathcal{Y}$-junction}
    \label{fig-y-junction}
\end{figure}

In order to verify this, we calculate the determinant of
$\mathbf{M}(\boldsymbol\lambda, \boldsymbol\beta)$ using the code in appendix. One such calculation is, for example, for $\lambda_{1,1}=0.44$, $\lambda_{1,2}=0.22$, $\beta_1=\beta_2=0.44 $ and, {by using the condition \eqref{hyp}}, $s$ varying between {$-0.06$ and $1.22$} (except $\frac{1}{2}$), we get $\det(\mathbf{M}(\boldsymbol\lambda, \boldsymbol\beta))=-3.2466 + 1.2854\ i$.

\begin{example}[Another unbalanced graph]

Now, we consider the Cauchy problem \eqref{KDV}-\eqref{condfront} on the case where
\[
a_1 = \dots = a_{k+m} = 1 \quad \text{and} \quad B = 
\begin{bmatrix} 
\frac{\sqrt{15}}{15} & \frac{\sqrt{15}}{15} & \frac{\sqrt{15}}{15} & \frac{\sqrt{15}}{15} & \frac{\sqrt{15}}{15} \\[10pt] 
\frac{\sqrt{15}}{15} & \frac{\sqrt{15}}{15} & \frac{\sqrt{15}}{15} & \frac{\sqrt{15}}{15} & \frac{\sqrt{15}}{15}\\[10pt]
\frac{\sqrt{15}}{15} & \frac{\sqrt{15}}{15} & \frac{\sqrt{15}}{15} & \frac{\sqrt{15}}{15} & \frac{\sqrt{15}}{15} 
\end{bmatrix}.\\[10pt]
\]
Here, we consider the star graph with three edges $(-\infty, 0)$ and five edges extending $(0, +\infty)$. Using Theorem \ref{grandeteorema} and analogous calculations, for each
$s \in \left(-0.06, 1.22 \right) \setminus \left\{\frac{1}{2}\right\}$, it is possible to determine $\lambda_{il}$ and $\beta_j$ such that the associated matrix $\mathbf{M}(\boldsymbol\lambda, \boldsymbol\beta)$ becomes invertible.
Thus, there exists a time $T > 0$ and a solution $(u_1, \dots, u_k, v_1, \dots, v_m)$  for the problem $\eqref{KDV}$, subject to the boundary conditions $\eqref{condfront}$ and satisfying the compatibility relation $u_{0,1}(0) = a_i u_{0,i}(0) = a_{k+j} v_{0,j}(0) \quad \text{for} \quad s > \frac{1}{2}$. As in the previous example, we found this through computer calculations with $\lambda_{1,1}=\lambda_{2,1}=\lambda_{3,1}=0.44$, $\lambda_{1,2}=\lambda_{2,2}=\lambda_{3,2}=0.22$ and $\beta_1=\beta_2=\beta_3=\beta_4=\beta_5=0.44$ (then $s\in(-0.06,1.22)\setminus\{\frac{1}{2}\}$), resulting in $det(\mathbf{M}(\boldsymbol\lambda, \boldsymbol\beta))= 6.2557 + 4.5450\ i$.  

\end{example}

\begin{figure}[h]
    \centering
\begin{tikzpicture}[line cap=round]
\def\numLeftEdges{5}
\def\numRightEdges{9}
\def\angleSpread{90}
\foreach \i in {1,...,\numLeftEdges} {
  \draw[thick, blue] (0,0) -- ({180 - \angleSpread/2 + (\i-1)*\angleSpread/(\numLeftEdges-1)}:3); 
  \draw[dashed, blue] ({180 - \angleSpread/2 + (\i-1)*\angleSpread/(\numLeftEdges-1)}:3) -- ({180 - \angleSpread/2 + (\i-1)*\angleSpread/(\numLeftEdges-1)}:4) node[left, blue] {};
  \draw[blue, <-, thick] ({180 - \angleSpread/2 + (\i-1)*\angleSpread/(\numLeftEdges-1)}:1.5) -- ({180 - \angleSpread/2 + (\i-1)*\angleSpread/(\numLeftEdges-1)}:2); 
  \draw[blue, <-, thick] ({180 - \angleSpread/2 + (\i-1)*\angleSpread/(\numLeftEdges-1)}:3.5) -- ({180 - \angleSpread/2 + (\i-1)*\angleSpread/(\numLeftEdges-1)}:3.75);
}
\foreach \i in {1,...,\numRightEdges} {
  \draw[thick, red] (0,0) -- ({\angleSpread/2 - (\i-1)*\angleSpread/(\numRightEdges-1)}:3); 
  \draw[dashed, red] ({\angleSpread/2 - (\i-1)*\angleSpread/(\numRightEdges-1)}:3) -- ({\angleSpread/2 - (\i-1)*\angleSpread/(\numRightEdges-1)}:4) node[right, red] {};
  \draw[red, ->, thick] ({\angleSpread/2 - (\i-1)*\angleSpread/(\numRightEdges-1)}:1.5) -- ({\angleSpread/2 - (\i-1)*\angleSpread/(\numRightEdges-1)}:2);
  \draw[red, ->, thick] ({\angleSpread/2 - (\i-1)*\angleSpread/(\numRightEdges-1)}:3.5) -- ({\angleSpread/2 - (\i-1)*\angleSpread/(\numRightEdges-1)}:3.75);
}
\filldraw (0,0) circle (2pt);
\end{tikzpicture}
\caption{A star graph with five $(-\infty,0)$ edges (the blue ones) and eight $(0,\infty)$ edges (the red ones)}
\end{figure}
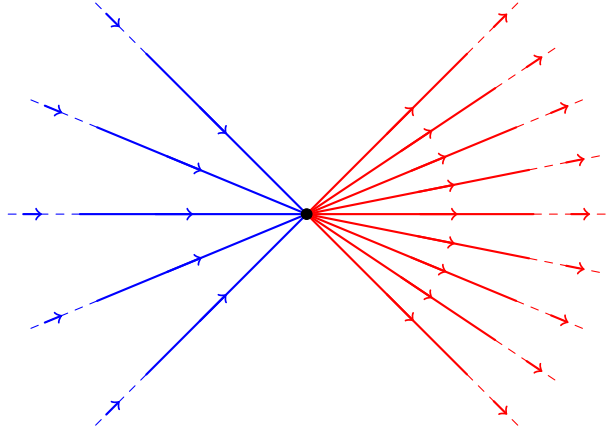

\begin{remark}
There are many other examples where the Cauchy problem \eqref{KDV}-\eqref{condfront} has solutions. For example, considering the parameters 
$a_1=\hdots=a_{k+m}=1$, $b_{ij}=\frac{\sqrt{km}}{km}$ $\forall i=1,\hdots,k\quad \forall j=1,\hdots,m$. Then choosing
$\lambda_{i1}=\beta_j=0.44$, $\lambda_{i2}=0.22$ - $s\in(-0.06,1.22)\setminus\{\frac{1}{2}\}$  -,  computational calculations can be performed to find $\mathbf{M}(\boldsymbol\lambda, \boldsymbol\beta)$ equal to:
\begin{itemize}
\item $-9.6406 + 1.2179\ i$ for $k=5$ and $m=9$;
\item $-3.2334 - 0.8302\ i$ for $k=11$ and $m=7$;
\item  $1.2747 + 0.9261\ i$ for $k=20$ and $m=30$.
\end{itemize}
\end{remark}

\begin{example}[Balanced graphs]\label{balanced}
Here, we consider the Cauchy problem \eqref{KDV}-\eqref{condfront} on the case that {$a_1=a_2=\hdots=a_{k+m}=1$ and $B= I_k$}. Let us consider the case of balanced star graphs. More precisely, we suppose $k=m$ to be an integer number between $1$ and $50$.  Assume {$s\in(-0.06,1.22)\setminus\{\frac{1}{2}\}$} a fixed real number and $u_{0,i}\in H^s(-\infty,0)$ and $v_{0,j}\in H^s(0,+\infty)$. 
 Then, there exists a time $T>0$ and a solution $(u_1,\hdots,u_k,v_1,\hdots,v_m)$ in the space $C([0,T];H^s(\mathcal{G}))$ (in the sense of distribution theory) of the problem $\eqref{KDV}$ with boundary conditions $\eqref{condfront}$ satisfying the compatibility conditions $u_{0,1}(0) = a_i u_{0,i}(0) = a_{k+j} v_{0,j}(0)$ when $s>\frac{1}{2}$. 
\end{example}

The claim in Example \ref{balanced} can be proved by calculating the determinant $\mathbf{M}(\boldsymbol\lambda, \boldsymbol\beta)$. We used the code in appendix to conclude that number is non-zero under corollary's hypothesis. Table \ref{tabeladassimulacoes} illustrates some of those calculations for $\lambda_{i1}=0.44$, $\lambda_{i2}=0.22$ and $\beta_j=0.44$. 

\begin{table}[H]
\begin{tabular}{ccccccllll}
 $k=m$ & $abs(det(\mathbf{M}(\boldsymbol\lambda, \boldsymbol\beta)))$ & $k=m$ & $abs(det(\mathbf{M}(\boldsymbol\lambda, \boldsymbol\beta)))$ & $k=m$ & $abs(det(\mathbf{M}(\boldsymbol\lambda, \boldsymbol\beta)))$\\
 1 & 2.37309 &     18 & 0.79410 &     35 & 0.02870 \cr
        2 & 3.75438 &     19 & 0.66305 &     36 & 0.02336 \cr
        3 & 4.45474 &     20 & 0.55210 &     37 & 0.01899 \cr
        4 & 4.69845 &     21 & 0.45856 &     38 & 0.01543 \cr
        5 & 4.64577 &     22 & 0.38001 &     39 & 0.01252 \cr
        6 & 4.40994 &     23 & 0.31426 &     40 & 0.01016 \cr
        7 & 4.06979 &     24 & 0.25940 &     41 & 0.00824 \cr
        8 & 3.67924 &     25 & 0.21374 &     42 & 0.00668 \cr
        9 & 3.27419 &     26 & 0.17584 &     43 & 0.00541 \cr
       10 & 2.87776 &     27 & 0.14445 &     44 & 0.00438 \cr
       11 & 2.50403 &     28 & 0.11849 &     45 & 0.00354 \cr
       12 & 2.16084 &     29 & 0.09708 &     46 & 0.00286 \cr
       13 & 1.85173 &     30 & 0.07944 &     47 & 0.00231 \cr
       14 & 1.57745 &     31 & 0.06493 &     48 & 0.00187 \cr
       15 & 1.33694 &     32 & 0.05302 &     49 & 0.00151 \cr
       16 & 1.12806 &     33 & 0.04325 &     50 & 0.00122 \cr
       17 & 0.94810 &     34 & 0.03525 &     51 & 0.00098
 
\end{tabular}
\caption{Computer Simulations for $\lambda_{i1}=0.44$, $\lambda_{i2}=0.22$ and $\beta_j=0.44$}
\label{tabeladassimulacoes}
\end{table}

\begin{remark}
In the examples above, one can expand the range in which $s$ can be taken. To do this, it is enough to look for different $\lambda_{il}$ and $\beta_j$ parameters for which the matrix $\mathbf{M}(\boldsymbol\lambda, \boldsymbol\beta)$ is invertible.
\end{remark}

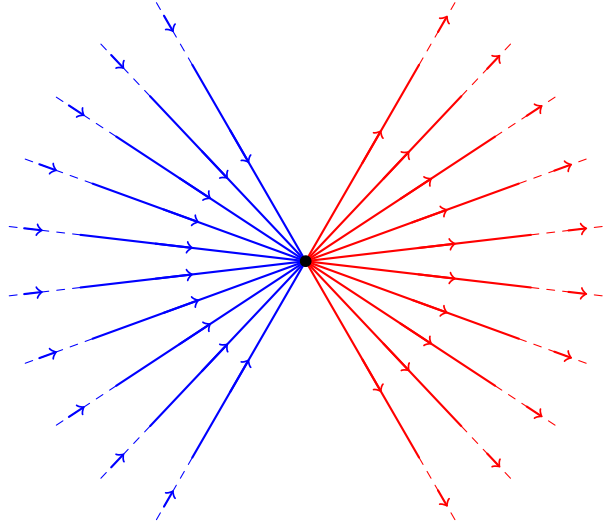
\begin{figure}[H]
    \centering
\begin{tikzpicture}[line cap=round]
\def\numLeftEdges{10}
\def\numRightEdges{10}
\def\angleSpread{120}
\foreach \i in {1,...,\numLeftEdges} {
  \draw[thick, blue] (0,0) -- ({180 - \angleSpread/2 + (\i-1)*\angleSpread/(\numLeftEdges-1)}:3); 
  \draw[dashed, blue] ({180 - \angleSpread/2 + (\i-1)*\angleSpread/(\numLeftEdges-1)}:3) -- ({180 - \angleSpread/2 + (\i-1)*\angleSpread/(\numLeftEdges-1)}:4) node[left, blue] {};
  \draw[blue, <-, thick] ({180 - \angleSpread/2 + (\i-1)*\angleSpread/(\numLeftEdges-1)}:1.5) -- ({180 - \angleSpread/2 + (\i-1)*\angleSpread/(\numLeftEdges-1)}:2); 
  \draw[blue, <-, thick] ({180 - \angleSpread/2 + (\i-1)*\angleSpread/(\numLeftEdges-1)}:3.5) -- ({180 - \angleSpread/2 + (\i-1)*\angleSpread/(\numLeftEdges-1)}:3.75);
}
\foreach \i in {1,...,\numRightEdges} {
  \draw[thick, red] (0,0) -- ({\angleSpread/2 - (\i-1)*\angleSpread/(\numRightEdges-1)}:3); 
  \draw[dashed, red] ({\angleSpread/2 - (\i-1)*\angleSpread/(\numRightEdges-1)}:3) -- ({\angleSpread/2 - (\i-1)*\angleSpread/(\numRightEdges-1)}:4) node[right, red] {};
  \draw[red, ->, thick] ({\angleSpread/2 - (\i-1)*\angleSpread/(\numRightEdges-1)}:1.5) -- ({\angleSpread/2 - (\i-1)*\angleSpread/(\numRightEdges-1)}:2);
  \draw[red, ->, thick] ({\angleSpread/2 - (\i-1)*\angleSpread/(\numRightEdges-1)}:3.5) -- ({\angleSpread/2 - (\i-1)*\angleSpread/(\numRightEdges-1)}:3.75);
}
\filldraw (0,0) circle (2pt);
\end{tikzpicture}
\caption{A balanced star graph with 20 edges}
\end{figure}

\section*{Acknowledgments}
 M. Cavalcante wishes to thank the support of CNPq, 
grant \# 310271/
2021-5, and the support of the Funda\c c\~ao de Amparo \`a Pesquisa do Estado de Alagoas
- FAPEAL, Brazil, grant \# E:60030.0000000161/2022. The authors thank the CAPES/Cofecub, grant \# 88887.879175/2023-00.
\section*{Statements and Declaration}

\textbf{Conflicts of interest} All authors declare that they have no conflicts of interest.

\appendix
\section{Code for computation of the determinant}
A MATLAB/Octave function to calculate the determinant of the matrix $\mathbf{M}(\boldsymbol\lambda, \boldsymbol\beta)$ follows. In input, $x$ means the $k$-dimensional vector $(\lambda_{1,1},\lambda_{2,1},\hdots,\lambda_{k1})$; $y$ means the $k$-dimensional vector $(\lambda_{1,2},\lambda_{2,2},\hdots,\lambda_{k2})$; $z$ means the $m$-dimensional vector $(\beta_{1},\beta_{2},\hdots,\beta_{m})$.\\

Furthermore, the meaning of $a(v),\hdots,h(v),j(v)$ is explained in the following table:

\begin{table}[H]
\begin{tabular}{ll}
 variable on code & variable on theorem\\
 $a(v)$ & $\rho_{v1}$\\[4pt]
 $b(v)$ & $\rho_{v2}$\\[4pt]
 $c(v)$ & $d_v$\\[4pt]
 $d(v)$ & $\tilde{\rho}_{v1}$\\[4pt]
 $e(v)$ & $\tilde{\rho}_{v2}$\\[4pt]
 $f(v)$ & $\tilde d_v$\\[4pt]
 $g(v)$ & $\tilde{\tilde{\rho}}_{v1}$\\[4pt]
 $h(v)$ & $\tilde{\tilde{\rho}}_{v2}$\\
 $j(v)$ & $\tilde{\tilde d}_v$\\
\end{tabular}
\end{table}

\begin{lstlisting}[language=MATLAB] 
% code for calculating the determinant
function retval = detergeral(x,y,z,A,B)

%number of edges
if isequal(size(x),size(y))
   k = max(size(y));
else
   disp('The vectors x and y are not the same size.')
end

m=max(size(z));

%building the matrix
for v = 1:m
	c(v)=exp(i*pi*z(v));
	f(v)=exp(i*pi*(z(v)-1));
	j(v)=exp(i*pi*(z(v)-2));
	end

for v = 1:k
	a(v)=2*sin(((pi)/(3))*x(v)+((pi)/(6)));
	d(v)=2*sin(((pi)/(3))*x(v)-((pi)/(6)));
	g(v)=2*sin(((pi)/(3))*x(v)-((pi)/(2)));
	b(v)=2*sin(((pi)/(3))*y(v)+((pi)/(6)));
	e(v)=2*sin(((pi)/(3))*y(v)-((pi)/(6)));
	h(v)=2*sin(((pi)/(3))*y(v)-((pi)/(2)));
	end

M=zeros(2*k+m);

% matrices M11 and M21 
M(1:k+m-1,1)=a(1)*ones(k+m-1,1);
M(1:k+m-1,2)=b(1)*ones(k+m-1,1);
for v = 2:k
	M(v-1,2*v-1) = -A(v)*a(v);
	end
for v = 2:k
	M(v-1,2*v) = -A(v)*b(v);
	end
% matrix M22
for v = 1:m
	M(k-1+v,2*k+v) = -A(k+v)*c(v);
	end
% matrix M31
for v = 1:k
	M(k+m-1+v,2*v-1) = -d(v);
	end
for v = 1:k
	M(k+m-1+v,2*v) = -e(v);
	end
% matrix M32
for v = 1:m
	M(k+m:2*k+m-1,2*k+v) = f(v)*B(1:k,v);
    end
% matrix M41
for v = 1:k
	M(2*k+m,2*v-1) = ((A(v))^(-1))*g(v);
	end
for v = 1:k
	M(2*k+m,2*v) = ((A(v))^(-1))*h(v);
    end
% matrix M42
for v = 1:m
	M(2*k+m,2*k+v) = -((A(k+v))^(-1))*j(v);
	end
retval=det(M)
end
\end{lstlisting}

\end{document}